\documentclass{amsart}
\usepackage{amssymb}
\usepackage{amsmath}
\usepackage{amsfonts}
\usepackage{amsthm}

\newtheorem{theorem}{Theorem}[section]
\newtheorem{lemma}[theorem]{Lemma}
\newtheorem{proposition}[theorem]{Proposition}
\newtheorem{corollary}[theorem]{Corollary}
\newtheorem{definition}[theorem]{Definition}
\newtheorem{example}[theorem]{Example}

\def\rmark{\mbox{$\rm\bf\rule{0.06em}{1.45ex}\kern-0.05em R$}}

\def\b1{b^{-1}}
\def\a1{a^{-1}}

\begin{document}

\begin{center}
{\textbf{\Large Study of Bipolar Fuzzy Soft Hypervector Spaces}}

\vspace{5mm}

{\textbf{O. R. Dehghan}}\\[0pt]
\vspace{0.3cm}

\begin{minipage}{12.5cm}
\vspace{0.5cm} \textbf{ Abstract.} {\small In this paper, three topics in bipolar fuzzy soft hypervector spaces are investigated. At first, four equivalent conditions to definition of a bipolar fuzzy soft hypervector space are presented, from different point of views. Then some new bipolar fuzzy soft hypervector spaces are constructed, using the notions of subhyperspace, level subset, generated subhyperspace and previous bipolar fuzzy soft hypervector spaces. Finally, normal bipolar fuzzy soft hypervector spaces is studied, as an especial kind of bipolar fuzzy soft hypervector spaces. All concepts are supported by interesting examples.}
\end{minipage}

\vspace{0.5cm}

\end{center}

\section{\protect\Large Introduction}
Fuzzy set was burn in 1965, when Zadeh \cite{Zadeh} defined it as a function $\nu:A\rightarrow [0,1]$. There are various generalizations for fuzzy sets. One of them is bipolar fuzzy set (Zhang \cite{Zhang}), which the image of membership function is increased from $[0,1]$ to $[-1,1]$.

Also, Molodtsov \cite{Molodtsov} introduced soft sets in 1999, that was an another mathematical concept for modeling vagueness. This idea was studied in different fields, especially in algebraic structures; for examples, soft groups, soft rings and soft vector spaces investigated by Aktas \cite{Aktas}, Acar \cite{Acar} and sezgin \cite{Sezgin}, respectively. Then, Cogman \cite{Cogman} defined a fuzzy soft set in 2011, that is a more accurate tool for modeling. After that, Abdollah \cite{Abdullah} introduced bipolar fuzzy soft sets in 2014 as an applied generalization of previous theories. In the following, some researchers applied this idea in different branches, for instance, Akram \cite{Akram,Akram 2} studied it in K-algebras, Ali \cite{Ali} used it in decision making, Abughazalah \cite{Abughazalah} applied the bipolar fuzzy sets in BCI-algebras, Riaz \cite{Riaz} discussed bipolar fuzzy soft topology, Mahmood \cite{Mahmood} presented it's complex extension and Khan \cite{Khan} studied bipolar fuzzy soft matrices.

On the other hand, Marty \cite{Marty} by generalization of operation into hyperoperation, introduced algebraic hyperstructures in 1934. An operation assigns to any two elements of $X$ a unique element of the set $X$, while a hyperoperation assigns them a unique subset of $X$. This idea has been studied in various fields, for more see the books \cite{Corsini 2}, \cite{Davvaz book} and \cite{Vougiouklis}. Particularly, Scafati-Tallini \cite{Tallini 1} introduced the notion of hypervector space in 1990. Hypervector spaces studied by Ameri \cite{Ameri Dehghan 1}, Sedghi \cite{Sedghi Dehghan} and the author \cite{Dehghan lin func hvs,Dehghan Neutro HVS,Dehghan Ameri Ebrahimi}.

The mentioned extensions of fuzzy sets have affected in algebraic hyperstructures (\cite{Davvaz book 2}). For example, fuzzy hypervector spaces introduced by Ameri \cite{Ameri FHVS over VF} in 2005. He investigated some properties of fuzzy hypervector spaces (\cite{Ameri Dehghan 2,Ameri Dehghan 3,Ameri Dehghan 4,Ameri Dehghan 5}). The author followed his work and studied some more results of fuzzy hypervector spaces (\cite{Dehghan Sum SP on FHVS,Dehghan Various FQHVS,Dehghan bal abs,Dehghan afin conv,Dehghan Norouzi points}). Ranjbar \cite{Ranjbar} checked out some properties of fuzzy soft hypervector spaces. The author \cite{Dehghan soft persian,Dehghan Nodehi soft} investigated some results in soft hypervector spaces. Norouzi \cite{Norouzi Ameri soft hmod} introduced some new directions on soft hypermodules and soft fuzzy hypermodules. Sarwar \cite{Sarwar} and Muhiuddin \cite{Muhiuddin} applied bipolar fuzzy soft sets in hypergraphs and hyper BCK-ideals, respectively.

Recently, the author \cite{Dehghan Int BF Soft HVS} applied bipolar fuzzy soft sets in hypervector spaces and studied some important results about bipolar fuzzy soft hypervector spaces. We follow \cite{Dehghan Int BF Soft HVS} in this paper and in section \ref{section equiv thm} present four equivalent conditions to defined notion. After that, we construct some new mentioned structure in section \ref{section new bf soft hvs}, using subhyperspaces, $(\alpha,\beta)$-level subsets and generated bipolar fuzzy soft sets. Next, we shortly study normal bipolar fuzzy soft hypervector spaces in section \ref{section normal bf soft hvs}. All new notions in the paper are supported by interesting examples. Finally, in section \ref{section conclution}, we present some ideas for next investigations.
\section{\protect\Large Preliminaries}
Here, some definitions and examples are presented from the published papers, for using in the rest of the article.
\begin{definition}\cite{Zhang}
If $\mathcal{X}$ is a non-empty set, $\mathcal{B}^{+}:\mathcal{X}\rightarrow [0,1]$ and $\mathcal{B}^{-}:\mathcal{X}\rightarrow [-1,0]$ show how much the member applies to the desired property and the implicit counter-property, respectively, then
\[\mathcal{B}=\{(x,\mathcal{B}^{+}(x),\mathcal{B}^{-}(x)),\ x\in \mathcal{X}\},\]
is called a bipolar fuzzy set (shortly, bf set) in $\mathcal{X}$. The collection of all bf sets over $\mathcal{U}$ is denoted by $BF^\mathcal{U}$.
\end{definition}
\begin{definition}\cite{Molodtsov}
If $\mathcal{B}$ is a set of parameters, $\mathcal{U}$ is the universe set with power set $P(\mathcal{U})$ and $\mathcal{G}:\mathcal{B}\rightarrow P(\mathcal{U})$ is a mapping, then a pair $(\mathcal{G},\mathcal{B})$ is called a soft set over $\mathcal{U}$.
\end{definition}
\begin{definition}\cite{Abdullah}
If $\mathcal{B}$ is a set of parameters, $\mathcal{U}$ is the universe set and $\mathcal{G}:\mathcal{B} \rightarrow BF^\mathcal{U}$ assigns to any parameter $e\in \mathcal{B}$ a bf set $\mathcal{G}_e$, i.e.
\[\forall e\in \mathcal{B};\ \mathcal{G}_{e}=\{(x,\mathcal{G}^{+}_{e}(x),\mathcal{G}^{-}_{e}(x)),\ x\in \mathcal{U}\},\]
then $(\mathcal{G},\mathcal{B})$ is said to be a bipolar fuzzy soft set (shortly, bfs set) of $\mathcal{U}$.
\end{definition}

\noindent Note that $(\mathcal{G},\mathcal{B})\sqsubseteq(\mathcal{H},\mathcal{C})$, iff $\mathcal{B}\subseteq \mathcal{C}$ and $\mathcal{G}_e^+(x)\leq \mathcal{G}_e^+(x)$, $\mathcal{G}_e^-(x)\geq \mathcal{G}_e^-(x)$, $\forall e\in \mathcal{B}, x\in \mathcal{U}$.
\begin{definition}\label{D HVS}\cite{Tallini 1}
A hypervector space (shortly, hvs) is an algebra $(\mathcal{V},+,\circ,\mathcal{K})$ that $(\mathcal{V},+)$ is an Abelian group, $\mathcal{K}$ is a field and $\circ:\mathcal{K}\times \mathcal{V}\rightarrow P_{\ast }(\mathcal{V})$ ($P_{\ast}(\mathcal{V})= P(\mathcal{V})\setminus\emptyset$) is an external hyperoperation (shortly, eho), such that
\begin{enumerate}
   \item[(H$_{1}$)] $b\circ (y+z)\subseteq b\circ y+b\circ z$,
   \item[(H$_{2}$)] $(b+c)\circ y\subseteq b\circ y+c\circ y$,
   \item[(H$_{3}$)] $b\circ (c\circ y)=(bc)\circ y$,
   \item[(H$_{4}$)] $b\circ (-y)=(-b)\circ y=-(b\circ y)$,
   \item[(H$_{5}$)] $y\in 1\circ y$.
\end{enumerate}
\noindent Note that in the right hand of (H$_{1}$), $b\circ y+b\circ z=\{r+s:r\in b\circ y, s\in b\circ z\}$. In a similar way, we have in (H$_{2}$). Moreover, $b\circ (c\circ y)=\underset{p\in c\circ y}{\bigcup}b\circ p$.

\noindent If equality holds in (H$_{1}$), then $\mathcal{V}$ is strongly right distributive (shortly, srd). A strongly left distributive (shortly, sld) hvs is similarly defined.
\end{definition}
If $\mathcal{V}_1,\mathcal{V}_2$ are hvs's over the field $\mathcal{K}$, such that $\mathcal{V}_1\subseteq \mathcal{V}_2$, then $\mathcal{V}_1$ is said to be a subhyperspace (shortly, shs) of $\mathcal{V}_2$, i.e. $y-z\in \mathcal{V}_1$ and $b\circ y\subseteq \mathcal{V}_1$, for all $b\in \mathcal{K}$, $y,z\in \mathcal{V}_1$.
\begin{example}\cite{Ameri Dehghan 3}\label{example hvs R3}
$(\mathbb{R}^{3},+,\circ,\mathbb{R})$ is a hvs, where $(\mathbb{R}^{3},+)$ is the real vector space and $a\circ(x_{0},y_{0},z_{0})=\{(ax_{0},ay_{0},z),\ z\in \mathbb{R}\}$.
\end{example}
\begin{example}\cite{Dehghan Int BF Soft HVS}\label{example hvs Z4}
Let $\mathcal{K}=\mathbb{Z}_2=\{0,1\}$ be the field of two numbers with the following operations:
\[
\begin{tabular}{ccc}
\begin{tabular}{c|c|c}
$+$ & $0$ & $1$ \\ \hline
$0$ & $0$ & $1$ \\ \hline
$1$ & $1$ & $0$ \\
\end{tabular}
& \ \ \ \ \ \ \ \ &
\begin{tabular}{c|c|c}
$\cdot $ & $0$ & $1$ \\ \hline
$0$ & $0$ & $0$ \\ \hline
$1$ & $0$ & $1$ \\
\end{tabular}
\end{tabular}
\]
Then $(\mathbb{Z}_4,+,\circ,\mathbb{Z}_2)$ is a hvs, such that the operation $``+:\mathbb{Z}_4\times \mathbb{Z}_4\rightarrow \mathbb{Z}_4"$ and the eho $``\circ:\mathbb{Z}_2\times \mathbb{Z}_4\rightarrow P_{*}(\mathbb{Z}_4)"$ are defined as follow:
\[
\begin{array}{ccc}
\begin{tabular}{c|c|c|c|c}
$+$ & $0$ & $1$ & $2$ & $3$ \\ \hline
$0$ & $0$ & $1$ & $2$ & $3$ \\ \hline
$1$ & $1$ & $2$ & $3$ & $0$ \\ \hline
$2$ & $2$ & $3$ & $0$ & $1$ \\ \hline
$3$ & $3$ & $0$ & $1$ & $2$ \\
\end{tabular}
&\ \ \ \ \  &
\begin{tabular}{c|c|c|c|c}
$\circ $ & $0$ & $1$ & $2$ & $3$ \\ \hline
$0$ & $\{0,2\}$ & $\{0\}$ & $\{0\}$ & $\{0\}$ \\ \hline
$1$ & $\{0,2\}$ & $\{1,2,3\}$ & $\{0,2\}$ & $\{1,2,3\}$
\end{tabular}
\end{array}
\]
\end{example}
\begin{definition}\cite{Dehghan Int BF Soft HVS}\label{def bf soft hvs}
Let $\mathcal{V}=(\mathcal{V},+,\circ,\mathcal{K})$ be a hvs and $(\mathcal{G},\mathcal{B})$ be a bfs set of $\mathcal{V}$. Then $(\mathcal{G},\mathcal{B})$ is a bipolar fuzzy soft hypervector space (shortly, bfs-hvs) of $\mathcal{V}$, iff
\begin{enumerate}
  \item $\mathcal{G}^{+}_e(y-z)\geq \mathcal{G}^{+}_e(y)\wedge \mathcal{G}^{+}_e(z)$, $\mathcal{G}^{-}_e(y-z) \leq \mathcal{G}^{-}_e(y)\vee \mathcal{G}^{-}_e(z)$,
  \item $\bigwedge\limits_{r\in b\circ y}\mathcal{G}^{+}_e(r)\geq \mathcal{G}^{+}_e(y)$, $\bigvee\limits_{r\in b\circ y}\mathcal{G}^{-}_e(r) \leq \mathcal{G}^{-}_e(y)$.
\end{enumerate}
\end{definition}
\begin{example}\cite{Dehghan Int BF Soft HVS}\label{example bf soft hvs R3}
$(\mathcal{G},\{a,b\})$ is a bfs-hvs of $\mathcal{V}=(\mathbb{R}^{3},+,\circ,\mathbb{R})$ defined in Example \ref{example hvs R3}, where $``\mathcal{G}^+_a,\mathcal{G}^+_b:\mathbb{R}^{3}\rightarrow [0,1]"$ and $``\mathcal{G}^-_a,\mathcal{G}^-_b:\mathbb{R}^{3}\rightarrow [-1,0]"$ are given by the followings:
\[
\mathcal{G}_{a}^{+}(t)=\left\{
\begin{array}{cl}
0.7 & t\in \{0\}\times \{0\}\times \mathbb{R}, \\
0.3 & t\in (\mathbb{R}\times \{0\}\times \mathbb{R})\setminus (\{0\}\times \{0\}\times \mathbb{R}), \\
0 & o.w.
\end{array}
\right.
\]

\[
\mathcal{G}_{a}^{-}(t)=\left\{
\begin{array}{cl}
-0.8 & t\in \{0\}\times \{0\}\times \mathbb{R}, \\
-0.4 & t\in (\mathbb{R}\times \{0\}\times \mathbb{R})\setminus (\{0\}\times \{0\}\times \mathbb{R}), \\
-0.2 & o.w.
\end{array}
\right.
\]

\[
\mathcal{G}_{b}^{+}(t)=\left\{
\begin{array}{cl}
0.9 & t\in \{0\}\times \{0\}\times \mathbb{R}, \\
0.4 & t\in (\mathbb{R}\times \{0\}\times \mathbb{R})\setminus (\{0\}\times \{0\}\times \mathbb{R}), \\
0.1 & o.w.
\end{array}
\right.
\]

\[
\mathcal{G}_{b}^{-}(t)=\left\{
\begin{array}{cl}
-0.6 & t\in \{0\}\times \{0\}\times \mathbb{R}, \\
-0.5 & t\in (\mathbb{R}\times \{0\}\times \mathbb{R})\setminus (\{0\}\times \{0\}\times \mathbb{R}), \\
-0.1 & o.w.
\end{array}
\right.
\]
\end{example}
\begin{example}\cite{Dehghan Int BF Soft HVS}\label{example bf soft hvs Z4}
$(\mathcal{G},\{c,d,e\})$ is a bfs-hvs of $\mathcal{V}=(\mathbb{Z}_4,+,\circ, \mathbb{Z}_2)$ defined in Example \ref{example hvs Z4}, where $``\mathcal{G}^+_c,\mathcal{G}^+_d,\mathcal{G}^+_e:\mathbb{Z}_4\rightarrow [0,1]"$ and $``\mathcal{G}^-_c,\mathcal{G}^-_d,\mathcal{G}^-_e:\mathbb{Z}_4\rightarrow [-1,0]"$ are given by the followings:
\[
\begin{array}{ccc}
\mathcal{G}_{c}^{+}(t)=\left\{
\begin{array}{cl}
0.5 & t=0,2, \\
0.3 & o.w.
\end{array}
\right.
& \ \ \ \ \  &
\mathcal{G}_{c}^{-}(t)=\left\{
\begin{array}{cl}
-0.4 & t=0,2, \\
-0.2 & o.w.
\end{array}
\right.
\end{array}
\]

\[
\begin{array}{ccc}
\mathcal{G}_{d}^{+}(t)=\left\{
\begin{array}{cl}
0.7 & t=0,2, \\
0.2 & o.w.
\end{array}
\right.
& \ \ \ \ \  &
\mathcal{G}_{d}^{-}(t)=\left\{
\begin{array}{cl}
-0.6 & t=0,2, \\
-0.3 & o.w.
\end{array}
\right.
\end{array}
\]

\[
\begin{array}{ccc}
\mathcal{G}_{e}^{+}(t)=\left\{
\begin{array}{cl}
0.8 & t=0,2, \\
0.4 & o.w.
\end{array}
\right.
& \ \ \ \ \  &
\mathcal{G}_{e}^{-}(t)=\left\{
\begin{array}{cl}
-0.7 & t=0,2, \\
-0.5 & o.w.
\end{array}
\right.
\end{array}
\]
\end{example}
\begin{definition}\cite{Dehghan Int BF Soft HVS}\label{def sum product bf soft set}
Let $(\mathcal{G},\mathcal{A})$ and $(\mathcal{H},\mathcal{B})$ be bfs sets of hvs $\mathcal{V}$ and $b\in \mathcal{K}$. Then the sum $(\mathcal{G},\mathcal{A})+(\mathcal{H},\mathcal{B})$ and the scalar product $b\circ(\mathcal{G},\mathcal{A})$ are defined as the bfs sets $(\mathcal{G}+\mathcal{H},\mathcal{A}\cap \mathcal{B})$ and $(a\circ \mathcal{G},\mathcal{A})$, respectively, where
\[(\mathcal{G}+\mathcal{H})^+_e(x)=\bigvee_{x=y+z}(\mathcal{G}^+_e(y)\wedge \mathcal{H}^+_e(z)),\]
\[(\mathcal{G}+\mathcal{H})^-_e(x)=\bigwedge_{x=y+z}(\mathcal{G}^-_e(y)\vee \mathcal{H}^-_e(z)),\]
\[
(b\circ \mathcal{G})_{e}^{+}(x)=
\left\{
\begin{array}{cl}
\bigvee\limits_{y\in b\circ r}\mathcal{G}_{e}^{+}(r) & \exists r\in \mathcal{V},y\in b\circ r, \\
0 & otherwise,
\end{array}
\right.
\]
\[
(b\circ \mathcal{G})_{e}^{-}(x)=
\left\{
\begin{array}{cl}
\bigwedge\limits_{y\in b\circ r}\mathcal{G}_{e}^{-}(r) & \exists r\in \mathcal{V},y\in b\circ r, \\
0 & otherwise.
\end{array}
\right.
\]
\end{definition}
\begin{lemma}\cite{Dehghan Int BF Soft HVS}\label{lem 2}
If $(\mathcal{G},\mathcal{B})$ and $(\mathcal{H},\mathcal{C})$ are bfs sets of hvs $\mathcal{V}=(\mathcal{V},+, \circ,\mathcal{K})$, then for all $e\in \mathcal{B}\cap \mathcal{C}$, $\acute{e}\in \mathcal{B}$, $y,z\in \mathcal{V}$, it follows that:

1) $(\mathcal{G}+\mathcal{H})_e^+(y+z)\geq \mathcal{G}_e^+(y)\wedge \mathcal{H}_e^+(z)$, $(\mathcal{G}+ \mathcal{H})_e^-(y+z)\leq \mathcal{G}_e^-(y)\vee \mathcal{H}_e^-(z)$.

2) $(\mathcal{G},\mathcal{B})\sqsubseteq 1\circ (\mathcal{G},\mathcal{B})$, $-(\mathcal{G},\mathcal{B}) \sqsubseteq (-1)\circ (\mathcal{G},\mathcal{B})$, where $(-\mathcal{G})_{\acute{e}}^+(y)= \mathcal{G}_{\acute{e}}^+(-y)$, $(-\mathcal{G})_{\acute{e}}^-(y)=\mathcal{G}_{\acute{e}}^-(-y)$.
\end{lemma}
\section{\protect\Large Equivalence Theorems}\label{section equiv thm}
In this section we present some equivalent conditions to definition of a bfs-hvs (Definition \ref{def bf soft hvs}). More precisely, these conditions are stated separately in theorems \ref{thm iff1}, \ref{thm iff2}, \ref{thm iff3}, \ref{thm iff4}, and are summarized in Corollary \ref{cor bf soft hvs}.
\begin{theorem}\label{thm iff1}
If $\mathcal{V}=(\mathcal{V},+,\circ,\mathcal{K})$ is a hvs and $(\mathcal{G},\mathcal{B})$ is a bfs set of $\mathcal{V}$, then $(\mathcal{G},\mathcal{B})$ is a bfs-hvs of $\mathcal{V}$  iff:
\begin{enumerate}
  \item $(\mathcal{G},\mathcal{B})+(\mathcal{G},\mathcal{B}))\sqsubseteq (\mathcal{G},\mathcal{B})$,
  \item $-(\mathcal{G},\mathcal{B})\sqsubseteq (\mathcal{G},\mathcal{B})$,
  \item $b\circ (\mathcal{G},\mathcal{B})\sqsubseteq (\mathcal{G},\mathcal{B})$, $\forall b\in \mathcal{K}$.
\end{enumerate}
\end{theorem}
\begin{proof}
Let $(\mathcal{G},\mathcal{B})$ be a bfs-hvs, $x\in \mathcal{V}$ and $e\in \mathcal{A}$. Then

\noindent 1) $(\mathcal{G}+\mathcal{G})_e^+(x)=\underset{x=y+z}{\bigvee}\mathcal{G}_e^+(y)\wedge \mathcal{G}_e^+(z)\leq \mathcal{G}_e^+(x)$, $(\mathcal{G}+\mathcal{G})_e^-(x)= \underset{x=y+z}{\bigwedge}\mathcal{G}_e^-(y)\vee \mathcal{G}_e^-(z)\geq \mathcal{G}_e^-(x)$.

\noindent 2) $(-\mathcal{G})_e^+(x)=\mathcal{G}_e^+(-x)\leq \mathcal{G}_e^+(-(-x))=\mathcal{G}_e^+(x)$ and $(-\mathcal{G})_e^-(x)=\mathcal{G}_e^-(-x)\geq \mathcal{G}_e^-(-(-x))=\mathcal{G}_e^-(x)$.

\noindent 3) If there does not exist $r\in \mathcal{V}$ with $x\in b\circ r$, then $(b\circ \mathcal{G})_e^+(x)= 0\leq \mathcal{G}_e^+(x)$ and $(b\circ \mathcal{G})_e^-(x)=0\geq \mathcal{G}_e^-(x)$. Also, if $r\in \mathcal{V}$ with $x\in b\circ r$, then $\mathcal{G}_e^+(r)\leq \underset{s\in b\circ r}{\bigwedge}\mathcal{G}_e^+(s)\leq \mathcal{G}_e^+(x)$ and $\mathcal{G}_e^-(r)\geq \underset{s\in b\circ r}{\bigvee} \mathcal{G}_e^-(s)\geq \mathcal{G}_e^-(x)$. Thus $(b\circ \mathcal{G})_e^+(x)=\underset{x\in b\circ r}{\bigvee} \mathcal{G}_e^+(t)\leq \mathcal{G}_e^+(x)$ and $(b\circ \mathcal{G})_e^-(x)=\underset{x\in b\circ r}{\bigwedge} \mathcal{G}_e^-(r)\geq \mathcal{G}_e^-(x)$.

Conversely, if $e\in \mathcal{B}$, $y,z\in \mathcal{V}$ and $b\in \mathcal{K}$, then

\noindent 1) By Lemma \ref{lem 2},
\begin{eqnarray*}
\mathcal{G}_{e}^{+}(y-z) &\geq &(\mathcal{G}+\mathcal{G})_{e}^{+}(y-z) \\
&\geq &\mathcal{G}_{e}^{+}(y)\wedge \mathcal{G}_{e}^{+}(-z) \\
&\geq &\mathcal{G}_{e}^{+}(y)\wedge (-\mathcal{G})_{e}^{+}(-z) \\
&=&\mathcal{G}_{e}^{+}(y)\wedge \mathcal{G}_{e}^{+}(z),
\end{eqnarray*}
and
\begin{eqnarray*}
\mathcal{G}_{e}^{-}(y-z) &\leq &(\mathcal{G}+\mathcal{G})_{e}^{-}(y-z) \\
&\leq &\mathcal{G}_{e}^{-}(y)\vee \mathcal{G}_{e}^{-}(-z) \\
&\leq &\mathcal{G}_{e}^{-}(y)\vee (-\mathcal{G})_{e}^{-}(-z) \\
&=&\mathcal{G}_{e}^{-}(y)\vee \mathcal{G}_{e}^{-}(z).
\end{eqnarray*}
\noindent 2) For all $r\in b\circ y$, $\mathcal{G}_{e}^{+}(r)\geq (b\circ \mathcal{G})_{e}^{+}(r)=\underset{r\in b\circ s}{\bigvee} \mathcal{G}_{e}^{+}(s)\geq \mathcal{G}_{e}^{+}(y)$ and $\mathcal{G}_{e}^{-}(r)\leq (b\circ \mathcal{G}) _{e}^{-}(r)=\underset{r\in b\circ s}{\bigwedge} \mathcal{G}_{e}^{-}(s)\leq \mathcal{G}_{e}^{-}(y)$. Hence $\underset{r\in b\circ y}{\bigwedge}\mathcal{G}_{e}^{+}(r)\geq \mathcal{G}_{e}^{+}(y)$ and $\underset{r\in b\circ y}{\bigwedge}\mathcal{G}_{e}^{-}(r)\leq \mathcal{G}_{e}^{-}(y)$.

Therefore, by Definition \ref{def bf soft hvs}, $(\mathcal{G},\mathcal{B})$ is a bfs-hvs of $\mathcal{V}$.
\end{proof}
\begin{lemma}\label{lem 1}
Suppose $\mathcal{V}=(\mathcal{V},+,\circ,\mathcal{K})$ is an invertible hvs (i.e. $x\in b\circ y\Rightarrow y\in b^{-1}\circ x$, $\forall b\in \mathcal{K}\setminus\{0\}$). If $(\mathcal{G},\mathcal{B})$ is a bfs-hvs of $\mathcal{V}$, then $(b\circ (\mathcal{G},\mathcal{B}))_e ^+ (y)\geq \mathcal{G}_e ^+ (y)$, $(b\circ (\mathcal{G},\mathcal{B}))_e ^- (y)\leq \mathcal{G}_e ^+ (y)$, for all $e\in \mathcal{B}$, $y\in \mathcal{V}$ and non-zero $b\in \mathcal{K}$.
\end{lemma}
\begin{proof}
By Definition \ref{D HVS}, $y\in 1\circ y=b\circ(b^{-1}\circ y)$, so $y\in b\circ r$ for some $r\in b^{-1}\circ y$. Then
\begin{equation*}
\left(b\circ(\mathcal{G},\mathcal{B})\right) _{e}^{+}(y)=\underset{y\in b\circ r}{\bigvee}\mathcal{G}_{e}^{+}(r)= \underset{r\in b^{-1}\circ y} {\bigvee}\mathcal{G}_{e}^{+}(r)\geq \underset{r\in b^{-1}\circ y}{\bigwedge} \mathcal{G}_{e}^{+}(r)\geq \mathcal{G}_{e}^{+}(y),
\end{equation*}
and
\begin{equation*}
\left(b\circ(\mathcal{G},\mathcal{B})\right)_{e}^{-}(y)=\underset{y\in b\circ r}{\bigwedge}\mathcal{G}_{e}^{-}(r) =\underset{r\in b^{-1}\circ y} {\bigwedge}\mathcal{G}_{e}^{-}(r)\leq \underset{r\in b^{-1}\circ y}{\bigvee} \mathcal{G}_{e}^{-}(r)\leq \mathcal{G}_{e}^{-}(y).
\end{equation*}
\end{proof}
\begin{proposition}
Assume $\mathcal{V}$ is invertible and srd-hvs. If $(\mathcal{G},\mathcal{C})$ and $(\mathcal{H},\mathcal{C})$ are bfs-hvs of $\mathcal{V}$, then for all non-zero $b\in \mathcal{K}$,
\[b\circ((\mathcal{G},\mathcal{C})+(\mathcal{H},\mathcal{C}))=b\circ(\mathcal{G},\mathcal{C})+b\circ(\mathcal{H}, \mathcal{C}).\]
\end{proposition}
\begin{proof}
Let $e\in\mathcal{C}$, $x\in \mathcal{V}$, $\alpha_{1}=\left(b\circ((\mathcal{G},\mathcal{C})+(\mathcal{H}, \mathcal{C}))\right)_{e}^{+}(x)$, $\alpha_{2}=\left(b\circ((\mathcal{G},\mathcal{C})+(\mathcal{H},\mathcal{C})) \right)_{e}^{-}(x)$,

\noindent $\beta_{1}=\left(b\circ(\mathcal{G},\mathcal{C})+ b\circ(\mathcal{G},\mathcal{C})\right) _{e}^{+}(x)$ and $\beta_{2}=\left(b\circ(\mathcal{G},\mathcal{C})+b\circ(\mathcal{H},\mathcal{C})\right) _{e}^{-}(x)$. More precisely,
\begin{equation*}
\alpha _{1}=\left\{
\begin{array}{cl}
\underset{x\in b\circ t}{\bigvee }(\mathcal{G}+\mathcal{H})_{e}^{+}(t) & \exists t\in \mathcal{V};x\in b\circ t, \\
0 & o.w.
\end{array}
\right.
\end{equation*}
\begin{equation*}
\alpha _{2}=\left\{
\begin{array}{cl}
\underset{x\in b\circ t}{\bigwedge }(\mathcal{G}+\mathcal{H})_{e}^{-}(t) & \exists t\in \mathcal{V};x\in b\circ t, \\
0 & o.w.
\end{array}
\right.
\end{equation*}
\begin{equation*}
\beta _{1}=\underset{x=y+z}{\bigvee }\left( b\circ (\mathcal{G},\mathcal{C})\right) _{e}^{+}(y)\wedge \left( b\circ(\mathcal{H},\mathcal{C})\right) _{e}^{+}(z),
\end{equation*}
\begin{equation*}
\beta _{2}=\underset{x=y+z}{\bigwedge }\left( b\circ (\mathcal{G},\mathcal{C})\right) _{e}^{-}(y)\vee \left( b\circ(\mathcal{H},\mathcal{C})\right) _{e}^{-}(z).
\end{equation*}
We must prove $\alpha_1=\alpha_2$ and $\beta_1=\beta_2$.

If $\nexists t\in \mathcal{V}$, with $x\in b\circ t$, then $\alpha_{1}=\alpha _{2}=0$. In this case, for any $y,z\in \mathcal{V}$, with $x=y+z$, either there does not exist $r\in \mathcal{V}$, with $y\in b\circ r$, or there does not exist $s\in \mathcal{V}$, with $z\in b\circ s$ (because if $r,s\in \mathcal{V}$ such that $y\in b\circ r$ and $z\in b\circ s$, then $x=y+z\in b\circ r+b\circ s=b\circ(r+s)$, a contradiction). Thus $\left( b\circ (\mathcal{G},\mathcal{C})\right)_{e}^{+}(y) =\left( b\circ (\mathcal{G},\mathcal{C})\right) _{e}^{-}(y)=0$ or $\left( b\circ(\mathcal{H},\mathcal{C})\right) _{e}^{+}(z)=\left( b\circ (\mathcal{H},\mathcal{C})\right) _{e}^{-}(z)=0$. So $\left( b\circ (\mathcal{G},\mathcal{C})\right) _{e}^{+}(y)\wedge \left( b\circ (\mathcal{H},\mathcal{C})\right) _{e}^{+}(z)=0$ and $\left( b\circ (\mathcal{G},\mathcal{C})\right) _{e}^{-}(y)\vee \left( b\circ (\mathcal{H},\mathcal{C})\right) _{e}^{-}(z)=0$. Hence, $\left( b\circ (\mathcal{G},\mathcal{C})+b\circ (\mathcal{H},\mathcal{C})\right) _{e}^{+}(x)=0$ and $\left( b\circ (\mathcal{G},\mathcal{C})+b\circ(\mathcal{H},\mathcal{C})\right) _{e}^{-}(x)=0$. Therefore $\beta _{1}=0=\alpha _{1}$ and $\beta_{2}=0=\alpha_{2}$.

Also, if $t\in \mathcal{V}$, with $x\in b\circ t$, then for any $\epsilon >0$, there exist $r,s\in \mathcal{V}$, $x\in b\circ r$, $x\in b\circ s$, with $(\mathcal{G}+\mathcal{H})_{e}^{+}(r)>\left( b\circ ((\mathcal{G},\mathcal{C})+ (\mathcal{H},\mathcal{C}))\right)_{e}^{+}(x)-\epsilon =\alpha _{1}-\epsilon $, and $(\mathcal{G}+\mathcal{H})_{e}^{-}(s)<\left(b\circ ((\mathcal{G},\mathcal{C})+(\mathcal{H},\mathcal{C}))\right) _{e}^{-}(x)+\epsilon =\alpha _{2}+\epsilon $. Thus there are $y_{1},y_{2},z_{1},z_{2}\in \mathcal{V}$, such that $y_{1}+z_{1}=r$, and $y_{2}+z_{2}=s$, such that $\mathcal{G}_{e}^{+}(y_{1})\wedge \mathcal{H}_{e}^{+}(z_{1})> \alpha_{1}-\epsilon $ and $\mathcal{G}_{e}^{-}(y_{2})\vee \mathcal{H}_{e}^{-}(z_{2})<\alpha_{2}+\epsilon $. Hence, $x\in b\circ r=b\circ(y_{1}+z_{1}) \subseteq b\circ y_{1}+b\circ z_{1}$, $x\in b\circ s=b\circ(y_{2}+z_{2})\subseteq b\circ y_{2}+b\circ z_{2}$, and so $x=\acute{y}_{1}+\acute{z}_{1}$, for some $\acute{y}_{1}\in b\circ y_{1}$, $\acute{z}_{1}\in b\circ z_{1}$, and $x=\acute{y}_{2}+\acute{z}_{2}$, for some $\acute{y} _{2}\in b\circ y_{2}$, $\acute{z}_{2}\in b\circ z_{2}$. By  Definition \ref{def bf soft hvs} and Lemma \ref{lem 1},
\begin{eqnarray*}
\beta _{1} &=&\left( b\circ (\mathcal{G},\mathcal{C})+b\circ (\mathcal{H},\mathcal{C})\right) _{e}^{+}(x) \\
&\geq &\left( b\circ (\mathcal{G},\mathcal{C})\right)_{e}^{+}(\acute{y}_{1})\wedge\left(b\circ(\mathcal{H}, \mathcal{C})\right) _{e}^{+}(\acute{z}_{1}) \\
&\geq &\mathcal{G}_{e}^{+}(\acute{y}_{1})\wedge \mathcal{H}_{e}^{+}(\acute{z}_{1}) \\
&\geq &\left( \underset{t_{1}\in b\circ y_{1}}{\bigwedge }\mathcal{G}_{e}^{+}(t_{1})\right) \wedge \left(\underset{l_{1}\in a\circ z_{1}}{\bigwedge }\mathcal{H}_{e}^{+}(l_{1})\right)  \\
&\geq &\mathcal{G}_{e}^{+}(y_{1})\wedge \mathcal{H}_{e}^{+}(z_{1}) \\
&>&\alpha _{1}-\epsilon,
\end{eqnarray*}
and
\begin{eqnarray*}
\beta _{2} &=&\left( b\circ (\mathcal{G},\mathcal{C})+b\circ (\mathcal{H},\mathcal{C})\right) _{e}^{-}(x) \\
&\leq &\left( b\circ(\mathcal{G},\mathcal{C})\right) _{e}^{-}(\acute{y}_{2})\vee \left( b\circ(\mathcal{H}, \mathcal{C}) \right) _{e}^{-}(\acute{z}_{2}) \\
&\leq &\mathcal{G}_{e}^{-}(\acute{y}_{2})\vee \mathcal{H}_{e}^{-}(\acute{z}_{2}) \\
&\leq &\left( \underset{t_{2}\in b\circ y_{2}}{\bigvee }\mathcal{G}_{e}^{-}(t_{2})\right) \vee \left( \underset{l_{2}\in a\circ z_{2}}{\bigvee }\mathcal{H}_{e}^{-}(l_{2})\right)  \\
&\leq &\mathcal{G}_{e}^{-}(y_{2})\vee \mathcal{G}_{e}^{-}(z_{2}) \\
&<&\alpha _{2}+\epsilon.
\end{eqnarray*}
Since $\epsilon >0$ was arbitrary, we have $\beta _{1}\geq \alpha _{1}$ and $\beta _{2}\leq \alpha _{2}$.

On the other hand, for all $\epsilon >0$, there exist $y_{1},y_{2},z_{1},z_{2}\in \mathcal{V}$, such that $x=y_{1}+z_{1} =y_{2}+z_{2}$, $\left( b\circ (\mathcal{G},\mathcal{C})\right) _{e}^{+}(y_{1})\wedge \left( b\circ (\mathcal{H},\mathcal{C})\right)_{e}^{+}(z_{1})>\beta _{1}-\epsilon $ and $\left( b\circ (\mathcal{G}, \mathcal{C}) \right)_{e}^{-} (y_{2}) \vee \left( b\circ (\mathcal{H},\mathcal{C})\right)_{e}^{-}(z_{2}) <\beta_{2}+\epsilon$. Taking $\beta _{1}>\epsilon $ and $\beta_{2}<\epsilon $ (if $\beta _{1}=0$, then $\alpha _{1}=0$ and if $\beta _{2}=0$, then $\alpha_{2}=0$, and there is nothing to prove), it follows that: $\left( b\circ(\mathcal{G},\mathcal{C})\right) _{e}^{+}(y_{1})>0$, $\left( b\circ (\mathcal{H},\mathcal{C}) \right)_{e}^{+}(z_{1})>0$, $\left( b\circ (\mathcal{G},\mathcal{C})\right) _{e}^{-}(y_{2})<0$ and $\left( b\circ (\mathcal{H},\mathcal{C})\right) _{e}^{-}(z_{2})<0$.  Thus there are $\acute{y}_{1}, \acute{z}_{1}, \acute{y}_{2},\acute{z}_{2}\in \mathcal{V}$, such that $y_{1}\in b\circ \acute{y}_{1}$, $z_{1}\in b\circ \acute{z}_{1}$, $y_{2}\in b\circ \acute{y}_{2}$ and $z_{2}\in b\circ \acute{z}_{2}$, with $\mathcal{G}_{e}^{+} (\acute{y} _{1})\wedge \mathcal{H}_{e}^{+} (\acute{z}_{1})>\beta _{1}-\epsilon $ and $\mathcal{G}_{e}^{-} (\acute{y}_{2})\vee \mathcal{H}_{e}^{-} (\acute{z}_{2}) <\beta _{2}+\epsilon $. But $x=y_{1}+z_{1}\in b\circ \acute{y}_{1}+b\circ \acute{z}_{1}=b\circ (\acute{y}_{1}+\acute{z}_{1})$, $x=y_{2}+z_{2}\in b\circ \acute{y}_{2}+b\circ \acute{z}_{2}=b\circ (\acute{y}_{2} +\acute{z}_{2})$, and so by Lemma \ref{lem 2}, $\alpha_{1}=\underset{x\in \text{$b\circ $}t}{\bigvee}(\mathcal{G}+\mathcal{H})_{e}^{+}(t)\geq(\mathcal{G}+ \mathcal{H})_{e}^{+}(\acute{y}_{1}+ \acute{z}_{1})\geq \mathcal{G}_{e}^{+} (\acute{y}_{1})\wedge \mathcal{H}_{e}^{+} (\acute{z}_{1})>\beta _{1}-\epsilon $, and $\alpha _{2}=\underset{x\in \text{$b\circ $}t}{\bigwedge }(\mathcal{G}+\mathcal{H})_{e}^{-}(t)\leq (\mathcal{G}+\mathcal{H})_{e}^{-} (\acute{y}_{2}+ \acute{z}_{2})\leq \mathcal{G}_{e}^{-}(\acute{y}_{2})\vee\mathcal{H}_{e}^{-}(\acute{z}_{2})<\beta_{2}+\epsilon$. Hence $\alpha _{1}\geq \beta_{1}$ and $\alpha _{2}\leq \beta _{2}$, since $\epsilon >0$ was arbitrary.

Therefore, the proof is completed.
\end{proof}
\begin{definition}\label{def level subset}
If $(\mathcal{G},\mathcal{B})$ is a bfs set of hvs $\mathcal{V}$, $\alpha\in(0,1]$, $\beta\in [-1,0)$, then the soft set
\[(\mathcal{G},\mathcal{B})_{\alpha,\beta}=\{\left(\mathcal{G}_e\right)_{\alpha,\beta};\ e\in \mathcal{B}\},\]
is called $(\alpha,\beta)$-level soft subset of $\mathcal{V}$, where
\[\left(\mathcal{G}_e\right)_{\alpha,\beta}=\{v\in \mathcal{V};\ \mathcal{G}_e^+(v)\geq \alpha,\ \mathcal{G}_e^-(v)\leq \beta\},\]
is an $(\alpha,\beta)$-level subset of the bfs set $\mathcal{G}_e=(\mathcal{G}_e^+,\mathcal{G}_e^-)$.
\end{definition}
\begin{example}
Let $(\mathcal{G},\mathcal{B})$ be the bfs set of the hvs $\mathcal{V}=(\mathbb{R}^{3},+,\circ,\mathbb{R})$ in Example \ref{example bf soft hvs R3}. Then
\begin{eqnarray*}
\left(\mathcal{G}_{a}\right)_{0.4,-0.3}&=&\{x\in\mathbb{R}^{3};\ \mathcal{G}_{a}^{+}(x)\geq 0.4,\ \mathcal{G}_{a}^{-}(x)\leq -0.3\} \\
&=&\left(\{0\}\times \{0\}\times \mathbb{R}\right) \cap \left( \mathbb{R}\times \{0\}\times \mathbb{R}\right)  \\
&=&\{0\}\times \{0\}\times \mathbb{R},
\end{eqnarray*}
and
\begin{eqnarray*}
\left(\mathcal{G}_{b}\right) _{0.4,-0.3} &=&\{x\in \mathbb{R}^{3};\ \mathcal{G}_{b}^{+}(x)\geq 0.4,\ \mathcal{G}_{b}^{-}(x)\leq -0.3\} \\
&=&\left(\mathbb{R}\times\{0\}\times\mathbb{R}\right)\cap \left(\mathbb{R}\times\{0\}\times \mathbb{R}\right) \\
&=&\mathbb{R}\times \{0\}\times \mathbb{R}.
\end{eqnarray*}
Thus
\begin{eqnarray*}
(\mathcal{G},\mathcal{B})_{0.4,-0.3} &=&\{\left(\mathcal{G}_{a}\right) _{0.4,-0.3},\left(\mathcal{G}_{b}\right) _{0.4,-0.3}\} \\
&=&\left\{\{0\}\times \{0\}\times \mathbb{R},\mathbb{R}\times \{0\}\times\mathbb{R}\right\}.
\end{eqnarray*}
Similarly,
\begin{equation*}
(\mathcal{G},\mathcal{B})_{0.8,-0.5}=\left\{\emptyset,\{0\}\times \{0\}\times \mathbb{R}\right\}.
\end{equation*}
\end{example}
\begin{example}
Let $\mathcal{V}=(\mathbb{R}^{3},+,\circ,\mathbb{R})$ be the hvs defined in Example \ref{example hvs R3}. Define a bfs set $(\mathcal{G},\mathcal{B})$ of $\mathcal{V}$, where $\mathcal{B}=\{a,b\}$, $``\mathcal{G}^+_a, \mathcal{G}^+_b:\mathbb{R}^{3}\rightarrow [0,1]"$ and $``\mathcal{G}^-_a,\mathcal{G}^-_b:\mathbb{R}^{3} \rightarrow [-1,0]"$ are given by the followings:
\begin{equation*}
\mathcal{G}_{a}^{+}(x,y,z)=\left\{
\begin{array}{cl}
0.3 & x\geq 0,yz\geq 0, \\
0.7 & x\geq 0,yz<0, \\
0.5 & x<0,
\end{array}
\right.
\end{equation*}
\begin{equation*}
\mathcal{G}_{a}^{-}(x,y,z)=\left\{
\begin{array}{cl}
-0.1 & y\geq 0,xz\geq 0, \\
-0.4 & y\geq 0,xz<0, \\
-0.6 & y<0,
\end{array}
\right.
\end{equation*}
\begin{equation*}
\mathcal{G}_{b}^{+}(x,y,z)=\left\{
\begin{array}{cl}
0.6 & x\geq 0,y\geq 0,z\geq 0, \\
0.7 & x\geq 0,y\geq 0,z<0, \\
0.3 & otherwise,
\end{array}
\right.
\end{equation*}
\begin{equation*}
\mathcal{G}_{b}^{-}(x,y,z)=\left\{
\begin{array}{cl}
-0.7 & z\geq 0, \\
-0.3 & z<0.
\end{array}
\right.
\end{equation*}
Then
\begin{equation*}
(\mathcal{G},\mathcal{B})_{0.6,-0.5}=\left\{\left(G_{a}\right)_{0.6,-0.5},\left(G_{b}\right)_{0.6,-0.5}\right\},
\end{equation*}
where,
\begin{equation*}
\left(\mathcal{G}_{a}\right)_{0.6,-0.5}=\{(x,y,z)\in \mathbb{R}^{3};\ x\geq 0,y<0,z>0\},
\end{equation*}
and
\begin{equation*}
\left(\mathcal{G}_{b}\right)_{0.6,-0.5}=\{(x,y,z)\in \mathbb{R}^{3};\ x\geq 0,y\geq 0,z\geq 0\}.
\end{equation*}
\end{example}
\begin{example}
Let $(\mathcal{G},\mathcal{C})$ be the bfs set of the hvs $\mathcal{V}=(\mathbb{Z}_{4},+,\circ,\mathbb{Z}_{2})$ in Example \ref{example bf soft hvs Z4}. Then
\begin{eqnarray*}
\left(\mathcal{G}_{c}\right)_{0.5,-0.5} &=&\{x\in \mathbb{Z}_{4};\ \mathcal{G}_{c}^{+}(x)\geq 0.5,\ \mathcal{G}_{c}^{-}(x)\leq -0.5\} \\
&=&\{0,2\} \cap \emptyset  \\
&=&\emptyset,
\end{eqnarray*}
\begin{eqnarray*}
\left(\mathcal{G}_{d}\right)_{0.5,-0.5} &=&\{x\in \mathbb{Z}_{4};\ \mathcal{G}_{d}^{+}(x)\geq 0.5,\ \mathcal{G}_{d}^{-}(x)\leq -0.5\} \\
&=&\{0,2\} \cap \{0,2\}  \\
&=&\{0,2\},
\end{eqnarray*}
and
\begin{eqnarray*}
\left(\mathcal{G}_{e}\right)_{0.5,-0.5} &=&\{x\in \mathbb{Z}_{4};\ \mathcal{G}_{e}^{+}(x)\geq 0.5,\ \mathcal{G}_{e}^{-}(x)\leq -0.5\} \\
&=&\{0,2\} \cap \{0,1,2,3\}  \\
&=&\{0,2\}.
\end{eqnarray*}
Thus
\begin{eqnarray*}
(\mathcal{G},\mathcal{C})_{0.5,-0.5} &=&\{\left(\mathcal{G}_{c}\right) _{0.5,-0.5},\left(\mathcal{G}_{d}\right) _{0.5,-0.5}, \left(\mathcal{G}_{e}\right) _{0.5,-0.5}\} \\
&=&\left\{\emptyset,\{0,2\},\{0,2\}\right\}.
\end{eqnarray*}
\end{example}
\begin{example}
Let $\mathcal{V}=(\mathbb{Z}_{4},+,\circ ,\mathbb{Z}_{2})$ be the hvs defined in Example \ref{example hvs Z4}. Consider a bfs set $(\mathcal{G},\mathcal{A})$ of $\mathcal{V}$, where $\mathcal{A}=\{c,d,e\}$ and $``\mathcal{G}_{c}^{+},\mathcal{G}_{d}^{+},\mathcal{G}_{e}^{+}: \mathbb{Z}_{4}\rightarrow [0,1]"$ and $``\mathcal{G}_{c}^{-},\mathcal{G}_{d}^{-},\mathcal{G}_{e}^{-}:\mathbb{Z}_{4}\rightarrow [-1,0]"$ are given by the followings:

\begin{equation*}
\begin{tabular}{lll}
\begin{tabular}{c|c|c|c|c}
$x$ & $0$ & $1$ & $2$ & $3$ \\ \hline
$\mathcal{G}_{c}^{+}(x)$ & $0.4$ & $0.3$ & $0.2$ & $0.7$
\end{tabular}
&  \ \ \ \  &
\begin{tabular}{c|c|c|c|c}
$x$ & $0$ & $1$ & $2$ & $3$ \\ \hline
$\mathcal{G}_{c}^{-}(x)$ & $-0.1$ & $-0.3$ & $-0.6$ & $-0.6$
\end{tabular}
\end{tabular}
\end{equation*}
\begin{equation*}
\begin{tabular}{lll}
\begin{tabular}{c|c|c|c|c}
$x$ & $0$ & $1$ & $2$ & $3$ \\ \hline
$\mathcal{G}_{d}^{+}(x)$ & $1$ & $0.5$ & $0.4$ & $0.1$
\end{tabular}
&  \ \ \ \  &
\begin{tabular}{c|c|c|c|c}
$x$ & $0$ & $1$ & $2$ & $3$ \\ \hline
$\mathcal{G}_{d}^{-}(x)$ & $-0.5$ & $-0.7$ & $-0.2$ & $-0.3$
\end{tabular}
\end{tabular}
\end{equation*}
\begin{equation*}
\begin{tabular}{lll}
\begin{tabular}{c|c|c|c|c}
$x$ & $0$ & $1$ & $2$ & $3$ \\ \hline
$\mathcal{G}_{e}^{+}(x)$ & $0.4$ & $0$ & $0.2$ & $0.8$
\end{tabular}
&  \ \ \ \  &
\begin{tabular}{c|c|c|c|c}
$x$ & $0$ & $1$ & $2$ & $3$ \\ \hline
$\mathcal{G}_{e}^{-}(x)$ & $0$ & $-0.1$ & $-0.7$ & $-0.5$
\end{tabular}
\end{tabular}
\end{equation*}
Then
\begin{eqnarray*}
\left(\mathcal{G}_{c}\right)_{0.3,-0.4} &=&\{x\in\mathbb{Z}_{4};\ \mathcal{G}_{c}^{+}(x)\geq 0.3,\mathcal{G} _{c}^{-}(x)\leq -0.4\} \\
&=&\{0,1,3\}\cap \{1,2,3\} \\
&=&\{1,3\},
\end{eqnarray*}
\begin{eqnarray*}
\left(\mathcal{G}_{d}\right)_{0.3,-0.4} &=&\{x\in \mathbb{Z}_{4};\ \mathcal{G}_{d}^{+}(x)\geq 0.3,\mathcal{G} _{d}^{-}(x)\leq -0.4\} \\
&=&\{0,1,2\}\cap \{0,1\} \\
&=&\{0,1\},
\end{eqnarray*}
and
\begin{eqnarray*}
\left(\mathcal{G}_{e}\right) _{0.3,-0.4} &=&\{x\in \mathbb{Z}_{4};\ \mathcal{G}_{e}^{+}(x)\geq 0.3,\mathcal{G} _{e}^{-}(x)\leq -0.4\} \\
&=&\{0,3\}\cap \{2,3\} \\
&=&\{3\}.
\end{eqnarray*}
Hence,
\begin{eqnarray*}
(\mathcal{G},\mathcal{A})_{0.3,-0.4} &=&\left\{ \left(\mathcal{G}_{c}\right)_{0.3,-0.4},\left(\mathcal{G} _{d}\right)_{0.3,-0.4},\left(\mathcal{G}_{e}\right) _{0.3,-0.4}\right\}  \\
&=&\left\{\{1,3\},\{0,1\},\{3\}\right\}.
\end{eqnarray*}
\end{example}
\begin{theorem}\label{thm iff2}
Consider $\mathcal{V}=(\mathcal{V},+,\circ,\mathcal{K})$ as a hypervector space and $(\mathcal{G},\mathcal{B})$ as a bfs set of $\mathcal{V}$. Then $(\mathcal{G},\mathcal{B})$ is classified as a bfs-hvs of $\mathcal{V}$ iff for all values of $\alpha$ within the range of $(0,1]$ and $\beta$ within the range of $[-1,0)$, the $(\alpha,\beta)$-level soft subset $(\mathcal{G},\mathcal{B})_{\alpha,\beta}$ qualifies as a soft hvs of $\mathcal{V}$, i.e. for all $\alpha\in(0,1]$, $\beta\in [-1,0)$ and $e\in \mathcal{B}$, $\left(\mathcal{G}_e\right)_{\alpha,\beta}$ is a shs of $\mathcal{V}$.
\end{theorem}
\begin{proof}
Suppose $(\mathcal{G},\mathcal{B})$ is a bfs-hvs of $\mathcal{V}$, $\alpha\in(0,1]$, $\beta\in [-1,0)$, $e\in \mathcal{B}$, $x,y\in \left(\mathcal{G}_e\right)_{\alpha,\beta}$ and $b\in \mathcal{K}$. Then $\mathcal{G}_e^+(x-y)\geq \mathcal{G}_e^+(x)\wedge \mathcal{G}_e^+(y)\geq \alpha$, $\mathcal{G}_e^-(x-y)\leq \mathcal{G}_e^-(x)\vee \mathcal{G}_e^-(y)\leq \beta$, and so $x-y\in \left(\mathcal{G}_e\right)_{\alpha,\beta}$. Also, for all $z\in b\circ x$,
\[\mathcal{G}_e^+(z)\geq \bigwedge\limits_{t\in b\circ x}\mathcal{G}_e^+(t)\geq \mathcal{G}_e^+(x)\geq \alpha,\]
and
\[\mathcal{G}_e^-(z)\leq \bigvee\limits_{t\in b\circ x}\mathcal{G}_e^-(t)\leq \mathcal{G}_e^-(x)\leq \beta.\]
Thus $z\in \left(\mathcal{G}_e\right)_{\alpha,\beta}$. Hence $b\circ x\subseteq \left(\mathcal{G}_e\right) _{\alpha,\beta}$. Therefore, $\left(\mathcal{G}_e\right)_{\alpha,\beta}$ is a shs of $\mathcal{V}$.

\noindent Conversely, suppose $\left(\mathcal{G}_e\right)_{\alpha,\beta}$ is a shs of $\mathcal{V}$, for all $\alpha\in(0,1]$, $\beta\in [-1,0)$ and $e\in \mathcal{A}$. Let $x,y\in \mathcal{V}$ and $b\in \mathcal{K}$. Choose $e\in \mathcal{B}$ and set $\alpha=\mathcal{G}_e^+(x)\wedge \mathcal{G}_e^+(y)$ and $\beta=\mathcal{G} _e^-(x)\vee \mathcal{G}_e^-(y)$. Then $\mathcal{G}_e^+(x),\mathcal{G}_e^+(y)\geq\alpha$, $\mathcal{G}_e^-(x), \mathcal{G}_e^-(y)\leq\beta$, and so $x,y\in\left(\mathcal{G}_e\right)_{\alpha,\beta}$. Thus $x-y\in \left(\mathcal{G}_e\right)_{\alpha,\beta}$. Hence $\mathcal{G}_e^+(x-y) \geq \alpha=\mathcal{G}_e^+(x)\wedge \mathcal{G}_e^+(y)$ and $\mathcal{G}_e^-(x-y) \leq \beta=\mathcal{G}_e^-(x)\vee \mathcal{G}_e^-(y)$.

\noindent Now choose $e\in \mathcal{B}$ and put $\alpha=\mathcal{G}_e^+(x)$ and $\beta=\mathcal{G}_e^-(x)$. Then $x\in\left(\mathcal{G}_e\right)_{\alpha,\beta}$ and so $b\circ x\subseteq \left(\mathcal{G}_e\right) _{\alpha,\beta}$. Thus
\[\bigwedge\limits_{t\in b\circ x}\mathcal{G}_e^+(t)\geq \bigwedge\limits_{s\in \left(\mathcal{G}_e\right) _{\alpha,\beta}} \mathcal{G}_e^+(s)\geq \alpha,\]
and
\[\bigvee\limits_{t\in b\circ x}\mathcal{G}_e^-(t)\leq \bigvee\limits_{s\in \left(\mathcal{G}_e\right) _{\alpha,\beta}} \mathcal{G}_e^-(s)\leq \beta.\]
Hence, $(\mathcal{G},\mathcal{B})$ is a bfs-hvs of $\mathcal{V}$.
\end{proof}
\begin{theorem}\label{thm iff3}
Suppose $(\mathcal{G},\mathcal{B})$ is a bfs set of sld-hvs $\mathcal{V}$. Then $(\mathcal{G},\mathcal{B})$ is a bfs-hvs of $\mathcal{V}$, iff
\[\bigwedge\limits_{t\in a\circ x+b\circ y}\mathcal{G}_e^+(t)\geq \mathcal{G}_e^+(x)\wedge \mathcal{G}_e^+(y),\]
and
\[\bigvee\limits_{t\in a\circ x+b\circ y}\mathcal{G}_e^-(t)\leq \mathcal{G}_e^-(x)\vee \mathcal{G}_e^-(y),\]
for all $e\in \mathcal{B}$, $x,y\in \mathcal{V}$, $a,b\in \mathcal{K}$.
\end{theorem}
\begin{proof}
If $(\mathcal{G},\mathcal{B})$ is a bfs-hvs of $\mathcal{V}$, $e\in \mathcal{B}$, $x,y\in \mathcal{V}$ and $a,b\in \mathcal{K}$, then
\begin{eqnarray*}
\bigwedge\limits_{t\in a\circ x+b\circ y}\mathcal{G}_e^+(t) &=& \bigwedge\limits_{t=t_1+t_2,t_1\in a\circ x,t_2\in b\circ y} \mathcal{G}_e^+(t) \\
&=& \bigwedge\limits_{t_1\in a\circ x,t_2\in b\circ y}\mathcal{G}_e^+(t_1+t_2) \\
&\geq&\left(\bigwedge\limits_{t_1\in a\circ x}\mathcal{G}_e^+(t_1)\right) \wedge \left(\bigwedge\limits_{t_2\in b\circ y} \mathcal{G}_e^+(t_2)\right) \\
&\geq& \mathcal{G}_e^+(x)\wedge \mathcal{G}_e^+(y),
\end{eqnarray*}
and
\begin{eqnarray*}
\bigvee\limits_{t\in a\circ x+b\circ y}\mathcal{G}_e^-(t) &=& \bigvee\limits_{t=t_1+t_2,t_1\in a\circ x,t_2\in b\circ y} \mathcal{G}_e^-(t) \\
&=& \bigvee\limits_{t_1\in a\circ x,t_2\in b\circ y}\mathcal{G}_e^-(t_1+t_2) \\
&\leq & \left(\bigvee\limits_{t_1\in a\circ x}\mathcal{G}_e^-(t_1)\right) \vee \left(\bigvee\limits_{t_2\in b\circ y} \mathcal{G}_e^- (t_2)\right) \\
&\leq & \mathcal{G}_e^-(x)\vee \mathcal{G}_e^-(y).
\end{eqnarray*}
Conversely, if $e\in \mathcal{B}$, $x,y\in \mathcal{V}$ and $a\in \mathcal{K}$, then by Definition \ref{D HVS},
\[\mathcal{G}_e^+(x-y)\geq \bigwedge\limits_{t\in 1\circ x-1\circ y}\mathcal{G}_e^+(t)\geq \mathcal{G}_e^+(x) \wedge \mathcal{G}_e^+(y),\]
and
\[G_e^-(x-y)\leq \bigvee\limits_{t\in 1\circ x-1\circ y}\mathcal{G}_e^-(t)\leq \mathcal{G}_e^-(x)\vee \mathcal{G} _e^-(y).\]
Also, $0\in 0\circ x$, for all $x\in \mathcal{V}$, since $\mathcal{V}$ is sld. Thus
\[\bigwedge\limits_{t\in a\circ x}\mathcal{G}_e^+(t)\geq \bigwedge\limits_{t\in 0\circ x+a\circ x}\mathcal{G} _e^+(t)\geq \mathcal{G}_e^+(x)\wedge \mathcal{G}_e^+(x)=\mathcal{G}_e^+(x),\]
and
\[\bigvee\limits_{t\in a\circ x}\mathcal{G}_e^-(t)\leq \bigvee\limits_{t\in 0\circ x+a\circ x}\mathcal{G} _e^-(t)\leq \mathcal{G}_e^-(x)\vee \mathcal{G}_e^-(x)=\mathcal{G}_e^-(x).\]
Therefore, by Definition \ref{def bf soft hvs}, $(\mathcal{G},\mathcal{B})$ is a bfs-hvs.
\end{proof}
\begin{theorem}\label{thm iff4}
Suppose $\mathcal{V}=(\mathcal{V},+,\circ,\mathcal{K})$ is a sld-hvs and $(\mathcal{G},\mathcal{B})$ is a bfs set of $\mathcal{V}$. Then $(\mathcal{G},\mathcal{B})$ is a bfs-hvs of $\mathcal{V}$ iff
\[\forall b,c\in \mathcal{K};\ b\circ (\mathcal{G},\mathcal{B})+c\circ (\mathcal{G},\mathcal{B}) \sqsubseteq (\mathcal{G},\mathcal{B}).\]
\end{theorem}
\begin{proof}
Let $(\mathcal{G},\mathcal{B})$ be a bfs-hvs of $\mathcal{V}$, $b,c\in \mathcal{K}$, $e\in \mathcal{B}$ and $x\in \mathcal{V}$. Then by Theorem \ref{thm iff1}, it follows that:
\begin{eqnarray*}
(b\circ \mathcal{G}+c\circ \mathcal{G})_e^+(x)&=&\underset{x=y+z}{\bigvee}((b\circ \mathcal{G})_e^+(y)\wedge (c\circ \mathcal{G})_e^+(z)) \\
&\leq &\underset{x=y+z}{\bigvee}(\mathcal{G}_e^+ (y)\wedge \mathcal{G}_e^+ (z)) \\
&=&(\mathcal{G}+\mathcal{G})_e^+(x) \\
&\leq &\mathcal{G}_e^+(x),
\end{eqnarray*}
\begin{eqnarray*}
(b\circ \mathcal{G}+c\circ \mathcal{G})_e^-(x)&=&\underset{x=y+z}{\bigwedge}((b\circ \mathcal{G})_e^-(y)\vee (c\circ \mathcal{G})_e^-(z)) \\
&\geq &\underset{x=y+z}{\bigwedge}(\mathcal{G}_e^- (y)\vee \mathcal{G}_e^- (z)) \\
&=&(\mathcal{G}+\mathcal{G})_e^-(x) \\
&\geq &\mathcal{G}_e^-(x).
\end{eqnarray*}
Thus $b\circ(\mathcal{G},\mathcal{B})+c\circ(\mathcal{G},\mathcal{B})\sqsubseteq (\mathcal{G},\mathcal{B})$.

\noindent Conversely, set $b=c=1$. Then by Lemma \ref{lem 2},
\[(\mathcal{G},\mathcal{B})+(\mathcal{G},\mathcal{B})\sqsubseteq 1\circ (\mathcal{G},\mathcal{B})+1\circ (\mathcal{G},\mathcal{B})\sqsubseteq (\mathcal{G},\mathcal{B}).\]
Now, put $b=1$ and $c=-1$. Then by Lemma \ref{lem 2}, for all $e\in \mathcal{B}$ and $x\in \mathcal{V}$, it follows that:
\begin{eqnarray*}
\mathcal{G}_e^+(0) &\geq & (1\circ (\mathcal{G},\mathcal{B})+(-1)\circ (\mathcal{G},\mathcal{B}))_e^+(0) \\
&=&\underset{0=y+z}{\bigvee}\left(\left(1\circ(\mathcal{G},\mathcal{B})\right)_e^+(y)\wedge\left((-1)\circ (\mathcal{G},\mathcal{B})\right)_e^+(z)\right)\\
&\geq &\left( 1\circ (\mathcal{G},\mathcal{B}) \right)_e^+ (x)\wedge \left( (-1)\circ (\mathcal{G},\mathcal{B}) \right)_e^+ (-x) \\
&\geq &\mathcal{G}_e^+ (x)\wedge (-\mathcal{G})_e^+(-x) \\
&=&\mathcal{G}_e^+ (x)\wedge \mathcal{G}_e^+(x) \\
&=&\mathcal{G}_e^+ (x),
\end{eqnarray*}
\begin{eqnarray*}
\mathcal{G}_e^-(0) &\leq & (1\circ (\mathcal{G},\mathcal{B})+(-1)\circ (\mathcal{G},\mathcal{B}))_e^-(0) \\
&=&\underset{0=y+z}{\bigwedge}\left(\left(1\circ(\mathcal{G},\mathcal{B})\right)_e^-(y)\vee\left((-1)\circ (\mathcal{G},\mathcal{B})\right)_e^-(z)\right)\\
&\leq &\left( 1\circ (\mathcal{G},\mathcal{B}) \right)_e^-(x)\vee \left((-1)\circ (\mathcal{G},\mathcal{B}) \right)_e^- (-x) \\
&\leq &\mathcal{G}_e^- (x)\vee (-\mathcal{G})_e^-(-x) \\
&=&\mathcal{G}_e^- (x)\vee \mathcal{G}_e^-(x) \\
&=&\mathcal{G}_e^- (x).
\end{eqnarray*}
So
\[
(0\circ\mathcal{G})_e^+(0)=\underset{0\in 0\circ t}{\bigvee}\mathcal{G}_e^+(t)=\mathcal{G}_e^+(0)=\underset{s\in \mathcal{V}}{\bigvee}\mathcal{G}_e^+(s),
\]
\[
(0\circ\mathcal{G})_e^-(0)=\underset{0\in 0\circ t}{\bigvee}\mathcal{G}_e^-(t)=\mathcal{G}_e^-(0)=\underset{s\in \mathcal{V}}{\bigwedge}\mathcal{G}_e^-(s).
\]
Now, suppose $b=0$ and $c=-1$. Then $0\circ (\mathcal{G},\mathcal{B}) +(-1)\circ (\mathcal{G},\mathcal{B}) \sqsubseteq (\mathcal{G},\mathcal{B})$, and so for all $e\in \mathcal{B}$, $x\in \mathcal{V}$, by Lemma \ref{lem 2}, we have:
\begin{eqnarray*}
\mathcal{G}_{e}^{+}(x) &\geq &\left(0\circ \mathcal{G}+(-1)\circ \mathcal{G}\right) _{e}^{+}(x) \\
&=&\underset{x=y+z}{\bigvee}\left(\left(0\circ\mathcal{G}\right)_{e}^{+}(y)\wedge\left((-1)\circ\mathcal{G} \right)_{e}^{+}(z)\right)  \\
&\geq &\left(0\circ \mathcal{G}\right)_{e}^{+}(0)\wedge\left((-1)\circ \mathcal{G}\right)_{e}^{+}(x) \\
&=&\mathcal{G}_{e}^{+}(0)\wedge \left((-1)\circ \mathcal{G}\right)_{e}^{+}(x) \\
&=&\left((-1)\circ \mathcal{G}\right)_{e}^{+}(x)\ \ \ \ \left(\underset{t\in\mathcal{V}}{\bigvee}\left((-1)\circ \mathcal{G}\right)_{e}^{+}(t)=\underset{t\in\mathcal{V}}{\bigvee}\mathcal{G}_{e}^{+}(t)=\mathcal{G}_{e}^{+}(0) \right) \\
&\geq &(-\mathcal{G})_{e}^{+}(x),
\end{eqnarray*}
\begin{eqnarray*}
\mathcal{G}_{e}^{-}(x) &\leq &\left(0\circ \mathcal{G}+(-1)\circ \mathcal{G}\right)_{e}^{-}(x) \\
&=&\underset{x=y+z}{\bigwedge}\left(\left(0\circ \mathcal{G}\right)_{e}^{-}(y)\vee\left((-1)\circ \mathcal{G} \right)_{e}^{-}(z)\right) \\
&\leq &\left(0\circ \mathcal{G}\right)_{e}^{-}(0)\vee \left((-1)\circ \mathcal{G}\right)_{e}^{-}(x) \\
&=&G_{e}^{-}(0)\vee \left((-1)\circ \mathcal{G}\right)_{e}^{-}(x) \\
&=&\left((-1)\circ \mathcal{G}\right)_{e}^{-}(x)\ \ \ \ \left(\underset{t\in\mathcal{V}}{\bigwedge}\left((-1) \circ \mathcal{G}\right)_{e}^{-}(t)=\underset{t\in \mathcal{V}}{\bigwedge}\mathcal{G}_{e}^{-}(t)=\mathcal{G} _{e}^{-}(0)\right) \\
&\leq &(-\mathcal{G})_{e}^{-}(x).
\end{eqnarray*}
Hence, $-(\mathcal{G},\mathcal{B})\sqsubseteq (\mathcal{G},\mathcal{B})$.

\noindent Now, by taking $c=0$, $b\circ (\mathcal{G},\mathcal{B})+0\circ(\mathcal{G},\mathcal{B}) \sqsubseteq (\mathcal{G},\mathcal{B})$. But for all $e\in \mathcal{B}$ and $x\in \mathcal{V}$,
\begin{eqnarray*}
((b\circ \mathcal{G})_e^++(0\circ \mathcal{G})_e^+)(x)&=&\underset{x=y+z}{\bigvee}((b\circ \mathcal{G})_e^+(y) \wedge (0\circ \mathcal{G})_e^+(z)) \\
&\geq&(b\circ \mathcal{G})_e^+(x)\wedge (0\circ \mathcal{G})_e^+(0) \\
&=&(b\circ \mathcal{G})_e^+(x)\wedge \mathcal{G}_e^+(0)\\
&=&(b\circ \mathcal{G})_e^+(x),\ \ \ \ \left(\underset{t\in \mathcal{V}}{\bigvee}(b\circ \mathcal{G})_e^+(t)= \underset{t\in \mathcal{V}}{\bigvee} \mathcal{G}_e^+(t)=\mathcal{G}_e^+(0)\right)
\end{eqnarray*}
\begin{eqnarray*}
((b\circ \mathcal{G})_e^-+(0\circ \mathcal{G})_e^-)(x)&=&\underset{x=y+z}{\bigwedge}((b\circ \mathcal{G})_e^-(y) \vee (0\circ \mathcal{G})_e^-(z)) \\
&\leq&(b\circ \mathcal{G})_e^-(x)\vee (0\circ \mathcal{G})_e^-(0) \\
&=&(b\circ \mathcal{G})_e^-(x)\vee \mathcal{G}_e^-(0)\\
&=&(b\circ \mathcal{G})_e^-(x).\ \ \ \ \left(\underset{t\in \mathcal{V}}{\bigwedge}(b\circ \mathcal{G})_e^-(t)= \underset{t\in \mathcal{V}}{\bigwedge} \mathcal{G}_e^-(t)=\mathcal{G}_e^-(0)\right)
\end{eqnarray*}
Thus $b\circ (\mathcal{G},\mathcal{B}) \sqsubseteq b\circ (\mathcal{G},\mathcal{B})+0\circ (\mathcal{G},\mathcal{B}) \sqsubseteq (\mathcal{G},\mathcal{B})$, for all $b \in \mathcal{K}$.

\noindent Therefore, by Theorem \ref{thm iff1}, $(\mathcal{G},\mathcal{B})$ is a bfs-hvs of $\mathcal{V}$.
\end{proof}
\begin{corollary}\label{cor bf soft hvs}
The following terms, about bfs set $(\mathcal{G},\mathcal{B})$ of $\mathcal{V}$ are equivalent:
\begin{enumerate}
  \item $(\mathcal{G},\mathcal{B})$ is a bfs-hvs of $\mathcal{V}$;
  \item $(\mathcal{G},\mathcal{B})+(\mathcal{G},\mathcal{B})\sqsubseteq (\mathcal{G},\mathcal{B})$, $-(\mathcal{G},\mathcal{B}) \sqsubseteq (\mathcal{G},\mathcal{B})$ and $b\circ (\mathcal{G},\mathcal{B}) \sqsubseteq (\mathcal{G},\mathcal{B})$, $\forall b\in \mathcal{K}$.
  \item For all $e\in \mathcal{B}$, $\alpha\in(0,1]$ and $\beta\in [-1,0)$, $\left(\mathcal{G}_e\right) _{\alpha,\beta}$ is a shs of $\mathcal{V}$;
  \end{enumerate}
Moreover, if $\mathcal{V}$ is sld, the above terms are equivalent to the followings:
\begin{enumerate}
  \item [(4)] $\bigwedge\limits_{r\in b\circ y+c\circ z}\mathcal{G}_e^+(r)\geq \mathcal{G}_e^+(y)\wedge \mathcal{G}_e^+(z)$, $\bigvee\limits _{r\in b\circ y+c\circ z}\mathcal{G}_e^-(r)\leq \mathcal{G}_e^-(y) \vee \mathcal{G}_e^-(z)$, $\forall e\in \mathcal{B}$, $y,z\in \mathcal{V}$, $b,c\in \mathcal{K}$.
  \item [(5)] $b\circ (\mathcal{G},\mathcal{B})+c\circ (\mathcal{G},\mathcal{B})\sqsubseteq (\mathcal{G},\mathcal{B})$, for all $b,c\in \mathcal{K}$.
\end{enumerate}
\end{corollary}
\begin{proof}
(1)$\Longleftrightarrow$(2): Theorem \ref{thm iff1};

(1)$\Longleftrightarrow$(3): Theorem \ref{thm iff2};

(1)$\Longleftrightarrow$(4): Theorem \ref{thm iff3};

(1)$\Longleftrightarrow$(5): Theorem \ref{thm iff4}.
\end{proof}
\section{\protect\Large New Constructed Bipolar Fuzzy Soft Hypervector Spaces}\label{section new bf soft hvs}
In this section, some ways to construct new bipolar fuzzy soft hypervector spaces are presented. In fact, in Proposition \ref{prop new bfshvs}, from every subhyperspace of $\mathcal{V}$, two bipolar fuzzy soft hypervector spaces are constructed. Also, in Theorem \ref{thm 1}, by using the notion of $(\alpha,\beta)$-level subsets, defined in Definition \ref{def level subset}, from every bfs-hvs $(\mathcal{F},\mathcal{A})$ of $\mathcal{V}$ and for every $e\in \mathcal{B}$, $\alpha\in (0,1]$, $\beta\in [-1,0)$, a new bfs-hvs of $\mathcal{V}$ is obtained. In particular, Theorem \ref{thm bf soft hvs generated}, provides a method to construct the bfs-hvs that is generated by a bfs set $(\mathcal{F},\mathcal{A})$ of $\mathcal{V}$.
\begin{proposition}\label{prop new bfshvs}
The following terms about a subset $\mathcal{X}$ of hvs $\mathcal{V}$ are equivalent:
\begin{enumerate}
  \item $\mathcal{X}$ is a shs of $\mathcal{V}$.
  \item $(\mathcal{G},\mathcal{B})$ is a bfs-hvs of $\mathcal{V}$, where $\mathcal{G}_e^+$, $\mathcal{G}_e^-$ are $\chi_\mathcal{X}$, the characteristic function of $\mathcal{X}$ and the zero function, respectively, for every $e\in \mathcal{B}$.
  \item $(\mathcal{H},\mathcal{B})$ is a bfs-hvs of $\mathcal{V}$, where
  \begin{equation*}
     \forall e\in \mathcal{B},\ \mathcal{H}_{e}^{+}(x)=0,\ \mathcal{H}_{e}^{-}(x)=\left\{
       \begin{array}{cc}
        -1 & x\in \mathcal{X}, \\
        0 & x\notin \mathcal{X}.
       \end{array}
     \right.
  \end{equation*}
\end{enumerate}
\end{proposition}
\begin{proof}
(1)$\Longleftrightarrow$(2): Let $\mathcal{X}$ be a shs of $\mathcal{V}$. Then for all $e\in \mathcal{B}$, $\alpha\in(0,1]$, $\beta\in[-1,0)$, $(\mathcal{G}_e)_{\alpha,\beta}$ is a shs of $\mathcal{V}$. Thus by Theorem \ref{thm iff2}, $(\mathcal{G},\mathcal{B})$ is a bfs-hvs of $\mathcal{V}$. Conversely, $(\mathcal{G}_e)_{1,\beta} =\mathcal{X}$ is a shs of $\mathcal{V}$, for arbitrary $e\in \mathcal{B}$ and $\beta\in[-1,0)$.

\noindent (1)$\Longleftrightarrow$(3): It is similar to (1)$\Leftrightarrow$(2).
\end{proof}
\begin{theorem}\label{thm 1}
If $(\mathcal{F},\mathcal{A})$ is a bfs-hvs of $\mathcal{V}$, then for any $\acute{e}\in \mathcal{A}$, $\alpha\in(0,1]$ and $\beta\in [-1,0)$, the bfs set $(\mathcal{G},\mathcal{A})$ of $\mathcal{V}$ given by
\begin{equation*}
\begin{array}{ccc}
\mathcal{G}_{e}^{+}(x)=\left\{
\begin{array}{cc}
\mathcal{F}_{e}^{+}(x) & x\notin \left(\mathcal{F}_{\acute{e}}\right) _{\alpha ,\beta } \\
1 & x\in \left(\mathcal{F}_{\acute{e}}\right) _{\alpha ,\beta }
\end{array}
\right.
& \text{ \ \ \ \ and \ \ \ \ } &
\mathcal{G}_{e}^{-}(x)=\left\{
\begin{array}{cc}
\mathcal{F}_{e}^{-}(x) & x\notin \left(\mathcal{F}_{\acute{e}}\right) _{\alpha ,\beta } \\
-1 & x\in \left(\mathcal{F}_{\acute{e}}\right) _{\alpha ,\beta }
\end{array}
\right.
\end{array}
\end{equation*}
is a bfs-hvs of $\mathcal{V}$.
\end{theorem}
\begin{proof}
Let's take a look at the requirements outlined in Definition \ref{def bf soft hvs} and evaluate if they are met. Suppose $\mathcal{W}=\left(\mathcal{F}_{\acute{e}}\right) _{\alpha ,\beta}$, $e\in \mathcal{A}$, $y,z\in \mathcal{V}$ and $a\in \mathcal{K}$.

\noindent 1) If $y,z\in \mathcal{W}$, then $y+z\in \mathcal{W}$ and so $\mathcal{G}_e^+(y+z)=1=\mathcal{G}_e^+(y) \wedge \mathcal{G}_e^+(z)$ and $\mathcal{G}_e^-(y+z)=-1=\mathcal{G}_e^-(y)\vee \mathcal{G}_e^-(z)$.

\noindent If $y\notin \mathcal{W}$, $z\notin \mathcal{W}$ and $y+z\in \mathcal{W}$, then $\mathcal{G}_e^+(y+z)= 1\geq \mathcal{G}_e^+(y)\wedge \mathcal{G}_e^+(z)$ and $\mathcal{G}_e^-(y+z)=-1\leq \mathcal{G}_e^-(y)\vee \mathcal{G}_e^-(z)$.

\noindent If $y\notin \mathcal{W}$, $z\notin \mathcal{W}$ and $y+z\notin \mathcal{W}$, then
\begin{eqnarray*}
\mathcal{G}_e^+(y+z)=\mathcal{F}_e^+(y+z)\geq\mathcal{F}_e^+(y)\wedge\mathcal{F}_e^+(z)=\mathcal{G}_e^+(y)\wedge \mathcal{G}_e^+(z),
\end{eqnarray*}
and
\begin{eqnarray*}
G_e^-(y+z)=\mathcal{F}_e^-(y+z)\leq \mathcal{F}_e^-(y)\vee \mathcal{F}_e^-(z)=\mathcal{G}_e^-(y)\vee \mathcal{G}_e^-(z).
\end{eqnarray*}

\noindent If $y\in \mathcal{W}$ and $z\notin \mathcal{W}$, then $y+z\notin \mathcal{W}$ (if $y+z=w\in \mathcal{W}$, then $z=w-y\in \mathcal{W}$, a contradiction), implies that
\begin{eqnarray*}
\mathcal{G}_{e}^{+}(y+z) &=&\mathcal{F}_{e}^{+}(y+z) \\
&\geq &\mathcal{F}_{e}^{+}(y)\wedge \mathcal{F}_{e}^{+}(z) \\
&\geq &\mathcal{F}_{e}^{+}(z) \\
&=&\mathcal{G}_{e}^{+}(z) \\
&=&1\wedge \mathcal{G}_{e}^{+}(z) \\
&=&\mathcal{G}_{e}^{+}(y)\wedge \mathcal{G}_{e}^{+}(z),
\end{eqnarray*}
and
\begin{eqnarray*}
\mathcal{G}_{e}^{-}(y+z) &=&\mathcal{F}_{e}^{-}(y+z) \\
&\leq &\mathcal{F}_{e}^{-}(y)\vee \mathcal{F}_{e}^{-}(z) \\
&\leq &\mathcal{F}_{e}^{-}(z) \\
&=&\mathcal{G}_{e}^{-}(z) \\
&=&-1\wedge \mathcal{G}_{e}^{-}(z) \\
&=&\mathcal{G}_{e}^{-}(y)\vee \mathcal{G}_{e}^{-}(z).
\end{eqnarray*}
\noindent If $y\notin \mathcal{W}$ and $z\in \mathcal{W}$, the result is similarly obtained.

\noindent Thus in any case, $\mathcal{G}_{e}^{+}(y+z)\geq \mathcal{G}_{e}^{+}(y)\wedge \mathcal{G}_{e}^{+}(z)$, and $\mathcal{G}_{e}^{-}(y+z)\leq \mathcal{G}_{e}^{-}(y) \vee \mathcal{G}_{e}^{-}(z)$. Also, by Theorem \ref{thm iff2}, $\mathcal{W}$ is a shs of $\mathcal{V}$ and so
\begin{equation*}
\mathcal{G}_{e}^{+}(-y)=\left\{
\begin{array}{cc}
\mathcal{F}_{e}^{+}(-y) & -x\notin \mathcal{W} \\
1 & -y\in \mathcal{W}
\end{array}
\right.
=\left\{
\begin{array}{cc}
\mathcal{F}_{e}^{+}(y) & x\notin \mathcal{W} \\
1 & y\in \mathcal{W}
\end{array}
\right.
=\mathcal{G}_{e}^{+}(y),
\end{equation*}
and
\begin{equation*}
\mathcal{G}_{e}^{-}(-y)=\left\{
\begin{array}{cc}
\mathcal{F}_{e}^{-}(-y) & -y\notin \mathcal{W} \\
-1 & -y\in \mathcal{W}
\end{array}
\right.
=\left\{
\begin{array}{cc}
\mathcal{F}_{e}^{-}(y) & x\notin \mathcal{W} \\
-1 & y\in \mathcal{W}
\end{array}
\right.
=\mathcal{G}_{e}^{-}(y).
\end{equation*}

\noindent 2) If $a\circ y\subseteq \mathcal{W}$, then $\bigwedge\limits_{t\in a\circ y}\mathcal{G}_{e}^{+}(t)= 1\geq \mathcal{G}_{e}^{+}(y)$  and $\bigvee \limits_{t\in a\circ y}\mathcal{G}_{e}^{-}(t)=-1\leq \mathcal{G} _{e}^{-}(y)$. If $a\circ y\nsubseteq \mathcal{W}$, then $y\notin \mathcal{W}$ (if $y\in \mathcal{W}$, then $a\circ y\subseteq \mathcal{W}$) and so
\begin{eqnarray*}
\bigwedge\limits_{t\in a\circ y}\mathcal{G}_{e}^{+}(t) &=&\left(\bigwedge\limits_{t\in a\circ y\cap \mathcal{W}} \mathcal{G}_{e}^{+}(t)\right) \wedge \left(\bigwedge\limits_{t\in a\circ y\setminus \mathcal{W}} \mathcal{G} _{e}^{+}(t)\right) \\
&=&\bigwedge\limits_{t\in a\circ y\setminus \mathcal{W}}\mathcal{G}_{e}^{+}(t) \\
&=&\bigwedge\limits_{t\in a\circ y\setminus \mathcal{W}}\mathcal{F}_{e}^{+}(t) \\
&\geq &\mathcal{F}_{e}^{+}(y) \\
&=&\mathcal{G}_{e}^{+}(y),
\end{eqnarray*}
and
\begin{eqnarray*}
\bigvee\limits_{t\in a\circ y}\mathcal{G}_{e}^{-}(t) &=&\left( \bigvee\limits_{t\in a\circ y\cap \mathcal{W}} \mathcal{G}_{e}^{-}(t)\right) \vee \left(\bigvee\limits_{t\in a\circ y\setminus \mathcal{W}}\mathcal{G} _{e}^{-}(t) \right) \\
&=&\bigvee\limits_{t\in a\circ y\setminus \mathcal{W}}\mathcal{G}_{e}^{-}(t) \\
&=&\bigvee\limits_{t\in a\circ y\setminus \mathcal{W}}\mathcal{F}_{e}^{-}(t) \\
&\leq &\mathcal{F}_{e}^{-}(y) \\
&=&\mathcal{G}_{e}^{-}(y).
\end{eqnarray*}
Therefore, $(\mathcal{G},\mathcal{A})$ is a bfs-hvs of $\mathcal{V}$.
\end{proof}
If $\mathcal{S}$ is a non-empty subset of a hypervector space $\mathcal{V}$, the linear span of $\mathcal{S}$ refers to the smallest shs of $\mathcal{V}$ that includes all the elements in $\mathcal{S}$. In other words, it is the intersection of all shs's $\mathcal{W}$ of $\mathcal{V}$ that contain $\mathcal{S}$. The linear span of $\mathcal{S}$ can be expressed as
\begin{equation*}
\left\langle \mathcal{S}\right\rangle =\left\{t\in \sum_{i=1}^{n}a{_{i}}\circ{s_{i}},a_{i}\in \mathcal{K},{s_{i}} \in \mathcal{S},n\in \mathbb{N}\right\}.
\end{equation*}
\begin{theorem}\label{thm bf soft hvs generated}
Suppose $(\mathcal{F},\mathcal{A})$ is a bfs set of hvs $\mathcal{V}=(\mathcal{V},+,\circ,\mathcal{K})$ such that $|Im (\mathcal{F}_e^+)|,|Im (\mathcal{F}_e^-)|<\infty$, for all $e\in \mathcal{A}$. Define subhyperspaces $\mathcal{W}_i$ of $\mathcal{V}$, $1\leq i\leq l$, by the followings:
\begin{equation*}
\mathcal{W}_{0}=\left\langle \mathcal{U}_{0}\right\rangle ,\ \ \mathcal{W}_{i}=\left\langle \mathcal{W}_{i-1}, \mathcal{U}_{i}\right\rangle,
\end{equation*}
where
\begin{equation*}
\mathcal{U}_{0}=\left\{x\in \mathcal{V},\ \mathcal{F}_{e}^{+}(x)=\bigvee\limits_{t\in \mathcal{V}}\mathcal{F} _{e}^{+}(t), \mathcal{F}_{e}^{-}(x)=\bigwedge\limits_{t\in \mathcal{V}}\mathcal{F}_{e}^{-}(t), \forall e\in \mathcal{A}\right\},
\end{equation*}
\begin{equation*}
\mathcal{U}_i=\left\{x\in \mathcal{V},\ \mathcal{F}_{e}^{+}(x)=\bigvee\limits_{t\in \mathcal{V}\setminus \mathcal{W}_{i-1}}\mathcal{F}_{e}^{+}(t), \mathcal{F}_{e}^{-}(x)=\bigwedge \limits_{t\in \mathcal{V}\setminus \mathcal{W}_{i-1}}\mathcal{F}_{e}^{-}(t),\forall e\in \mathcal{A}\right\},
\end{equation*}
$l\leq \left\vert Im(\mathcal{F}_{e}^{+})\right\vert,\left\vert Im(\mathcal{F}_{e}^{-})\right\vert$ and $\mathcal{W}_{l}=\mathcal{V}$. Then the bfs set $(\mathcal{G},\mathcal{A})$ of $\mathcal{V}$ defined by
\begin{equation*}
\left[
\begin{array}{cl}
\mathcal{G}: & \mathcal{A}\longrightarrow BF^{\mathcal{V}} \\
& \begin{array}{cl}
e\mapsto \mathcal{G}_{e} & :\mathcal{V}\rightarrow [0,1]\times [-1,0] \\
& \ \ x\mapsto (\mathcal{G}_{e}^{+}(x),\mathcal{G}_{e}^{-}(x)),
\end{array}
\end{array}
\right.
\end{equation*}
is the smallest bfs-hvs of $\mathcal{V}$ that contains $(\mathcal{F},\mathcal{A})$, where
\begin{equation*}
\mathcal{G}_{e}^{+}(x)=\left\{
\begin{array}{cl}
\bigvee\limits_{t\in \mathcal{V}}\mathcal{F}_{e}^{+}(t) & x\in \hat{\mathcal{W}}_{0}=\mathcal{W}_{0}, \\
\bigvee\limits_{t\in \mathcal{V}\setminus \mathcal{W}_{i-1}}\mathcal{F}_{e}^{+}(t) & x\in \hat{\mathcal{W}}_i= \mathcal{W}_{i}\setminus \mathcal{W}_{i-1},
\end{array}
\right.
\end{equation*}
and
\begin{equation*}
\mathcal{G}_{e}^{-}(x)=\left\{
\begin{array}{cl}
\bigwedge\limits_{t\in \mathcal{V}}\mathcal{F}_{e}^{-}(t) & x\in \hat{\mathcal{W}}_{0}=\mathcal{W}_{0}, \\
\bigwedge\limits_{t\in \mathcal{V}\setminus \mathcal{W}_{i-1}}\mathcal{F}_{e}^{-}(t) & x\in \hat{\mathcal{W}}_i= \mathcal{W}_{i}\setminus \mathcal{W}_{i-1}.
\end{array}
\right.
\end{equation*}
\end{theorem}
\begin{proof}
For all $\alpha\in(0,1]$, $\beta\in[-1,0)$ and $e\in \mathcal{A}$, $(\mathcal{G}_e)_{\alpha,\beta}$ is a subhyperspace of $\mathcal{V}$, so by Theorem \ref{thm iff2}, $(\mathcal{G},\mathcal{A})$ is a bfs-hvs of $\mathcal{V}$.

\noindent Now let $(\mathcal{H},\mathcal{A})$ be a bfs-hvs of $\mathcal{V}$ containing $(\mathcal{F}, \mathcal{A})$. We prove that $(\mathcal{F},\mathcal{A})\sqsubseteq (\mathcal{G},\mathcal{A})\sqsubseteq (\mathcal{H},\mathcal{A})$, i.e.
\[\forall e\in \mathcal{A}, \forall x\in \mathcal{V},\ \mathcal{G}_e^+(x)\leq \mathcal{H}_e^+(x),\ \mathcal{G} _e^-(x)\geq \mathcal{H}_e^-(x).\]
If $e\in \mathcal{A}$, $x\in \mathcal{W}_0$, then there exist $a_{01},a_{02},\ldots,a_{0n}\in \mathcal{K}$ and $t_{01},t_{02},\ldots,t_{0n}\in \mathcal{U}_0$ such that $x\in a_{01}\circ t_{01}+a_{02}\circ t_{02}+\cdots+ a_{0n}\circ t_{0n}$. Thus by Theorem \ref{thm iff3},
\begin{eqnarray*}
\mathcal{H}_{e}^{+}(x) &\geq &\bigwedge\limits_{t\in a_{01}\circ t_{01}+\cdots+a_{0n}\circ t_{0n}}\mathcal{H} _{e}^{+}(t) \\
&\geq &\bigwedge\limits_{t\in a_{01}\circ t_{01}+\cdots +a_{0n}\circ t_{0n}}\mathcal{F}_{e}^{+}(t) \\
&\geq &\mathcal{F}_{e}^{+}(t_{01})\wedge \cdots \wedge \mathcal{F}_{e}^{+}(t_{0n}) \\
&=&\bigvee\limits_{t\in \mathcal{V}}\mathcal{F}_{e}^{+}(t) \\
&\geq &\mathcal{F}_{e}^{+}(x),
\end{eqnarray*}
and
\begin{eqnarray*}
\mathcal{H}_{e}^{-}(x) &\leq &\bigvee\limits_{t\in a_{01}\circ t_{01}+\cdots+a_{0n}\circ t_{0n}}\mathcal{H} _{e}^{-}(t) \\
&\leq &\bigvee\limits_{t\in a_{01}\circ t_{01}+\cdots +a_{0n}\circ t_{0n}}\mathcal{F}_{e}^{-}(t) \\
&\leq &\mathcal{F}_{e}^{-}(t_{01})\vee \cdots \vee \mathcal{F}_{e}^{-}(t_{0n}) \\
&=&\bigwedge\limits_{t\in \mathcal{V}}\mathcal{F}_{e}^{-}(t) \\
&\leq &\mathcal{F}_{e}^{-}(x).
\end{eqnarray*}
Hence $\mathcal{F}_{e}^{+}(x)\leq \mathcal{G}_{e}^{+}(x)\leq \mathcal{H}_{e}^{+}(x)$ and $\mathcal{F}_{e}^{-}(x) \geq \mathcal{G}_{e}^{-}(x)\geq \mathcal{H}_{e}^{-}(x)$, for all $x\in \mathcal{W}_0$, $e\in \mathcal{A}$.

\noindent If $e\in \mathcal{A}$ and $x\in \hat{\mathcal{W}}_i$, $i=1,\ldots,l$, then there exist $a_{i1},a_{i2}, \ldots,a_{in}\in \mathcal{K}$ and $t_{i1},t_{i2},\ldots,t_{in}\in \mathcal{U}_i$ such that $x\in a_{i1}\circ t_{i1}+a_{i2}\circ t_{i2}+\cdots+a_{in}\circ t_{in}$, and so
\begin{eqnarray*}
\mathcal{H}_{e}^{+}(x) &\geq &\bigwedge\limits_{t\in a_{i1}\circ t_{i1}+\cdots+a_{in}\circ t_{in}}\mathcal{H} _{e}^{+}(t) \\
&\geq &\bigwedge\limits_{t\in a_{i1}\circ t_{i1}+\cdots +a_{in}\circ t_{in}}\mathcal{F}_{e}^{+}(t) \\
&\geq &\mathcal{F}_{e}^{+}(t_{i1})\wedge \cdots \wedge \mathcal{F}_{e}^{+}(t_{in}) \\
&=&\bigvee\limits_{t\in \mathcal{V}\setminus \mathcal{W}_{i-1}}\mathcal{F}_{e}^{+}(t) \\
&\geq &\bigvee\limits_{t\in \hat{\mathcal{W}}_{i}}\mathcal{F}_{e}^{+}(t) \\
&\geq &\mathcal{F}_{e}^{+}(x),
\end{eqnarray*}
and
\begin{eqnarray*}
\mathcal{H}_{e}^{-}(x) &\leq &\bigvee\limits_{t\in a_{i1}\circ t_{i1}+\cdots+a_{in}\circ t_{in}}\mathcal{H} _{e}^{-}(t) \\
&\leq &\bigvee\limits_{t\in a_{i1}\circ t_{i1}+\cdots +a_{in}\circ t_{in}}\mathcal{F}_{e}^{-}(t) \\
&\leq &\mathcal{F}_{e}^{-}(t_{i1})\vee \cdots \vee \mathcal{F}_{e}^{-}(t_{in}) \\
&=&\bigwedge\limits_{t\in \mathcal{V}\setminus \mathcal{W}_{i-1}}\mathcal{F}_{e}^{-}(t) \\
&\leq &\bigwedge\limits_{t\in \hat{\mathcal{W}}_{i}}\mathcal{F}_{e}^{-}(t) \\
&\leq &\mathcal{F}_{e}^{-}(x).
\end{eqnarray*}
Then $\mathcal{F}_{e}^{+}(x)\leq \mathcal{G}_{e}^{+}(x)\leq \mathcal{H}_{e}^{+}(x)$ and $\mathcal{F}_{e}^{-}(x) \geq \mathcal{G}_{e}^{-}(x)\geq \mathcal{H}_{e}^{-}(x)$, for all $x\in \hat{\mathcal{W}}_i$, $e\in \mathcal{A}$.

\noindent Therefore, $(\mathcal{G},\mathcal{A})$ is the smallest bfs-hvs of $\mathcal{V}$ containing $(\mathcal{F},\mathcal{A})$.
\end{proof}
\noindent $(\mathcal{G},\mathcal{A})$ defined in Theorem \ref{thm bf soft hvs generated}, is said to be the bipolar fuzzy soft hypervector space generated by $(\mathcal{F},\mathcal{A})$.
\section{\protect\Large Normal Bipolar Fuzzy Soft Hypervector Spaces}\label{section normal bf soft hvs}
Here a particular type of bfs-hvs's is shortly studied, supported by some examples. It is shown that we can obtain normal bfs-hvs's from bfs-hvs's, also from subhyperspaces of $\mathcal{V}$.
\begin{definition}
A bfs-hvs $(\mathcal{G},\mathcal{B})$ of $\mathcal{V}$ is considered normal if there exists $x$ in $\mathcal{V}$ with $\mathcal{G}_e^+(x)=1$ and $\mathcal{G}_e^-(x)=-1$, for all $e\in \mathcal{B}$.
\end{definition}
The fact is clearly evident, if $(\mathcal{G},\mathcal{B})$ is a bfs-hvs of $\mathcal{V}=(\mathcal{V},+,\circ, \mathcal{K})$, then $\mathcal{G}_e^+(x)\leq \mathcal{G}_e^+(0)$ and $\mathcal{G}_e^-(x)\geq \mathcal{G}_e^-(0)$, $\forall e\in \mathcal{B}$, $x\in \mathcal{V}$. Thus if $(\mathcal{G},\mathcal{B})$ is a normal bfs-hvs of $\mathcal{V}$,$\mathcal{G}_e^+(0)=1$ and $\mathcal{G}_e^-(0)=-1$, for all $e\in\mathcal{B}$. Hence $(\mathcal{G}, \mathcal{B})$ is a normal bfs-hvs of $\mathcal{V}$, iff $\mathcal{G}_e^+(0)=1$ and $\mathcal{G}_e^-(0)=-1$, for all $e\in \mathcal{B}$.
\begin{example}\label{example normal bf soft hvs R3}
Let $\mathcal{V}=(\mathbb{R}^{3},+,\circ,\mathbb{R})$ be the hvs in Example \ref{example hvs R3}, and $\mathcal{A}=\{a,b\}$ be the parameters. Define $``\mathcal{F}^+_a,\mathcal{F}^+_b:\mathbb{R}^{3}\rightarrow [0,1]"$ and $``\mathcal{F}^-_a,\mathcal{F}^-_b:\mathbb{R}^{3}\rightarrow [-1,0]"$ by the followings:
\[
\mathcal{F}_{a}^{+}(t)=\left\{
\begin{array}{cl}
1 & t\in \{0\}\times \{0\}\times \mathbb{R}, \\
0.4 & t\in (\mathbb{R}\times \{0\}\times \mathbb{R})\setminus (\{0\}\times \{0\}\times \mathbb{R}), \\
0 & otherwise,
\end{array}
\right.
\]

\[
\mathcal{F}_{a}^{-}(t)=\left\{
\begin{array}{cl}
-1 & t\in \{0\}\times \{0\}\times \mathbb{R}, \\
-0.6 & t\in (\mathbb{R}\times \{0\}\times \mathbb{R})\setminus (\{0\}\times \{0\}\times \mathbb{R}), \\
-0.2 & otherwise,
\end{array}
\right.
\]

\[
\mathcal{F}_{b}^{+}(t)=\left\{
\begin{array}{cl}
1 & t\in \{0\}\times \{0\}\times \mathbb{R}, \\
0.5 & t\in (\mathbb{R}\times \{0\}\times \mathbb{R})\setminus (\{0\}\times \{0\}\times \mathbb{R}), \\
0.3 & otherwise,
\end{array}
\right.
\]

\[
\mathcal{F}_{b}^{-}(t)=\left\{
\begin{array}{cl}
-1 & t\in \{0\}\times \{0\}\times \mathbb{R}, \\
-0.9 & t\in (\mathbb{R}\times \{0\}\times \mathbb{R})\setminus (\{0\}\times \{0\}\times \mathbb{R}), \\
-0.1 & otherwise.
\end{array}
\right.
\]
Then $(\mathcal{F},\mathcal{A})$ is a normal bfs-hvs of $\mathbb{R}^{3}$.
\end{example}
\begin{example}\label{example normal bf soft hvs Z4}
Let $\mathcal{V}=(\mathbb{Z}_4,+,\circ,\mathbb{Z}_2)$ be the hvs in Example \ref{example hvs Z4}, and $\mathcal{A}=\{c,d,e\}$ be the parameters. Define $``\mathcal{F}^+_c,\mathcal{F}^+_d,\mathcal{F}^+_e: \mathbb{Z}_4\rightarrow [0,1]"$ and $``\mathcal{F}^-_c,\mathcal{F}^-_d,\mathcal{F}^-_e:\mathbb{Z}_4\rightarrow [-1,0]"$ by:
\[
\begin{array}{ccc}
\mathcal{F}_{c}^{+}(t)=\left\{
\begin{array}{cc}
1 & t=0,2, \\
0.4 & t=1,3,
\end{array}
\right.
& \ \ \ \ \  &
\mathcal{F}_{c}^{-}(t)=\left\{
\begin{array}{cc}
-1 & t=0,2, \\
-0.3 & t=1,3,
\end{array}
\right.
\end{array}
\]

\[
\begin{array}{ccc}
\mathcal{F}_{d}^{+}(t)=\left\{
\begin{array}{cc}
1 & t=0,2, \\
0.6 & t=1,3,
\end{array}
\right.
& \ \ \ \ \  &
\mathcal{F}_{d}^{-}(t)=\left\{
\begin{array}{cc}
-1 & t=0,2, \\
-0.2 & t=1,3,
\end{array}
\right.
\end{array}
\]

\[
\begin{array}{ccc}
\mathcal{F}_{e}^{+}(t)=\left\{
\begin{array}{cc}
1 & t=0,2, \\
0.1 & t=1,3,
\end{array}
\right.
& \ \ \ \ \  &
\mathcal{F}_{e}^{-}(t)=\left\{
\begin{array}{cc}
-1 & t=0,2, \\
-0.8 & t=1,3.
\end{array}
\right.
\end{array}
\]
Then $(\mathcal{F},\mathcal{A})$ is a normal bfs-hvs of $\mathbb{Z}_4$.
\end{example}
\begin{proposition}
If $(\mathcal{G},\mathcal{B})$ is a bfs-hvs of $\mathcal{V}=(\mathcal{V},+,\circ,\mathcal{K})$, then the bfs set $(\tilde{\mathcal{G}},\mathcal{B})$ of $\mathcal{V}$ is a normal bfs-hvs of $\mathcal{V}$, containing $(\mathcal{G},\mathcal{B})$, where
\[\forall e\in \mathcal{B},\ \tilde{\mathcal{G}}_e^+(x)=\mathcal{G}_e^+(x)+1-\mathcal{G}_e^+(0),\ \tilde{\mathcal{G}}_e^-(x)=\mathcal{G}_e^-(x)-1+\mathcal{G}_e^-(0).\]
\end{proposition}
\begin{proof}
Let $e\in \mathcal{B}$, $y,z\in \mathcal{V}$ and $b\in \mathcal{K}$. Then
\begin{eqnarray*}
  \tilde{\mathcal{G}}_e^+(y-z) &=& \mathcal{G}_e^+(y-z)+1-\mathcal{G}_e^+(0) \\
    &\geq& (\mathcal{G}_e^+(y)\wedge \mathcal{G}_e^+(z))+1-\mathcal{G}_e^+(0) \\
    &=& (\mathcal{G}_e^+(y)+1-\mathcal{G}_e^+(0))\wedge (\mathcal{G}_e^+(z)+1-\mathcal{G}_e^+(0)) \\
    &=& \tilde{\mathcal{G}}_e^+(y)\wedge \tilde{\mathcal{G}}_e^+(z),
\end{eqnarray*}
and
\begin{eqnarray*}
  \tilde{\mathcal{G}}_e^-(y-z) &=& \mathcal{G}_e^-(y-z)-1+\mathcal{G}_e^-(0) \\
    &\leq& (\mathcal{G}_e^-(y)\vee \mathcal{G}_e^-(z))-1+\mathcal{G}_e^-(0) \\
    &=& (\mathcal{G}_e^-(y)-1+\mathcal{G}_e^-(0))\vee (\mathcal{G}_e^-(z)-1+\mathcal{G}_e^-(0)) \\
    &=& \tilde{\mathcal{G}}_e^-(y)\vee \tilde{\mathcal{G}}_e^-(z).
\end{eqnarray*}
Also,
\[\bigwedge\limits_{t\in b\circ y}\tilde{\mathcal{G}}_e^+(t)=\bigwedge\limits_{t\in b\circ y}(\mathcal{G}_e^+(t)+ 1-\mathcal{G}_e^+(0)) \geq \mathcal{G}_e^+(y)+1-\mathcal{G}_e^+(0)=\tilde{\mathcal{G}}_e^+(y),\]
and
\[\bigvee\limits_{t\in b\circ y}\tilde{\mathcal{G}}_e^-(t)=\bigvee\limits_{t\in b\circ y}(\mathcal{G}_e^-(t)-1+ \mathcal{G}_e^-(0)) \leq \mathcal{G}_e^-(y)-1+\mathcal{G}_e^-(0)=\tilde{\mathcal{G}}_e^-(y).\]
Hence, by Definition \ref{def bf soft hvs}, $(\tilde{\mathcal{G}},\mathcal{B})$ is a bfs-hvs of $\mathcal{V}$. Moreover, $\tilde{\mathcal{G}}_e^+(0)=1$ and $\tilde{\mathcal{G}}_e^-(0)=-1$, for all $e\in \mathcal{B}$. Next, for all $e\in \mathcal{B}$ and $x\in \mathcal{V}$ we have, $\tilde{\mathcal{G}}_e^+(y)=\mathcal{G}_e^+(y)+1- \mathcal{G}_e^+(0)\geq \mathcal{G}_e^+(y)$ and $\tilde{\mathcal{G}}_e^-(y)=\mathcal{G}_e^-(y)-1+\mathcal{G} _e^-(0) \leq \mathcal{G}_e^-(y)$. Thus $(\mathcal{G},\mathcal{B})\sqsubseteq (\tilde{\mathcal{G}},\mathcal{B})$. The proof is now completed.
\end{proof}
Clearly, a bfs-hvs $(\mathcal{G},\mathcal{B})$ of $\mathcal{V}$ is normal, iff $(\tilde{\mathcal{G}}, \mathcal{B})=(\mathcal{G},\mathcal{B})$.
\begin{proposition}
If $\mathcal{X}$ is a shs of $\mathcal{V}$, then the bfs set $(\mathcal{G},\mathcal{B})$ of $\mathcal{V}$ given by the following is a normal bfs-hvs of $\mathcal{V}$.
\begin{equation*}
\forall e\in \mathcal{B},\ \mathcal{G}_e^+(x)=
\left\{
\begin{array}{cc}
  1 & x\in \mathcal{X}, \\
  0 & x\notin \mathcal{X},
\end{array}
\right.
\ \ \
\text{and}
\ \ \
\mathcal{G}_e^-(x)=
\left\{
\begin{array}{cc}
  -1 & x\in \mathcal{X}, \\
  0 & x\notin \mathcal{X}.
\end{array}
\right.
\end{equation*}
\end{proposition}
\begin{proof}
Clearly, for all $e\in \mathcal{B}$, $x\in \mathcal{X}$, $\mathcal{G}_e^+(x)=1$ and $\mathcal{G}_e^-(x)=-1$.
\end{proof}
\begin{theorem}
If $(\mathcal{G},\mathcal{B})$ is a non-zero bfs-hvs of $V$, then $(\breve{\mathcal{G}},B)$ given by the following is a normal bfs-hvs of $V$, containing $(\mathcal{G},\mathcal{B})$.
\[\forall e\in B, x\in V,\ \breve{\mathcal{G}}_e^+(x)=\frac{\mathcal{G}_e^+(x)}{\mathcal{G}_e^+(0)},\ \breve{\mathcal{G}}_e^-(x)=-\frac{\mathcal{G}_e^-(x)} {\mathcal{G}_e^-(0)}.\]
\end{theorem}
\begin{proof}
Let $e\in \mathcal{B}$, $y,z\in \mathcal{V}$ and $b\in \mathcal{K}$. Then
\begin{eqnarray*}
\breve{\mathcal{G}}_{e}^{+}(y-z) &=&\frac{\mathcal{G}_{e}^{+}(y-z)}{\mathcal{G}_{e}^{+}(0)} \\
&\geq &\frac{\mathcal{G}_{e}^{+}(y)\wedge \mathcal{G}_{e}^{+}(z)}{\mathcal{G}_{e}^{+}(0)} \\
&=&\frac{\mathcal{G}_{e}^{+}(y)}{\mathcal{G}_{e}^{+}(0)}\wedge \frac{\mathcal{G}_{e}^{+}(z)} {\mathcal{G}_{e}^{+}(0)} \\
&=&\breve{\mathcal{G}}_{e}^{+}(y)\wedge \breve{\mathcal{G}}_{e}^{+}(z),
\end{eqnarray*}
and
\begin{eqnarray*}
\breve{\mathcal{G}}_{e}^{-}(y-z) &=&-\frac{\mathcal{G}_{e}^{-}(y-z)}{\mathcal{G}_{e}^{-}(0)} \\
&\leq &-\frac{\mathcal{G}_{e}^{-}(y)\vee \mathcal{G}_{e}^{-}(z)}{\mathcal{G}_{e}^{-}(0)} \\
&=&\left(-\frac{\mathcal{G}_{e}^{-}(y)}{\mathcal{G}_{e}^{-}(0)}\right) \vee \left(-\frac{\mathcal{G}_{e}^{-}(z)} {\mathcal{G}_{e}^{-}(0)}\right) \\
&=&\breve{\mathcal{G}}_{e}^{-}(y)\vee \breve{\mathcal{G}}_{e}^{-}(z).
\end{eqnarray*}
Also,
\begin{eqnarray*}
\bigwedge\limits_{t\in b\circ y}\breve{\mathcal{G}}_{e}^{+}(t)&=&\bigwedge\limits_{t\in b\circ y} \frac{\mathcal{G}_{e}^{+}(t)} {\mathcal{G}_{e}^{+}(0)} \\
&=&\frac{\bigwedge\limits_{t\in b\circ y}\breve{\mathcal{G}}_{e}^{+}(t)}{\mathcal{G}_{e}^{+}(0)}\\
&\geq &\frac{\mathcal{G}_{e}^{+}(y)}{\mathcal{G}_{e}^{+}(0)} \\
&=&\breve{\mathcal{G}}_{e}^{+}(y),
\end{eqnarray*}
and
\begin{eqnarray*}
\bigvee\limits_{t\in b\circ y}\breve{\mathcal{G}}_{e}^{-}(t) &=&\bigvee\limits_{t\in b\circ y} \left(-\frac{\mathcal{G}_{e}^{-}(t)} {\mathcal{G}_{e}^{-}(0)}\right)  \\
&=&-\frac{\bigvee\limits_{t\in b\circ y}\breve{\mathcal{G}}_{e}^{-}(t)}{\mathcal{G}_{e}^{-}(0)}\\
&\leq &-\frac{\mathcal{G}_{e}^{-}(y)}{\mathcal{G}_{e}^{-}(0)} \\
&=&\breve{\mathcal{G}}_{e}^{-}(y).
\end{eqnarray*}
Thus by Definition \ref{def bf soft hvs}, $(\breve{\mathcal{G}},\mathcal{B})$ is a bfs-hvs of $\mathcal{V}$. Moreover, $\breve{\mathcal{G}}_{e}^{+}(0)=1$, $\breve{\mathcal{G}}_{e}^{-}(0)=-1$, and for all $y\in \mathcal{V}$, $\mathcal{G}_{e}^{+}(y)\leq \breve{\mathcal{G}}_{e}^{+}(y)$, $\mathcal{G}_{e}^{-}(y)\geq \breve{\mathcal{G}}_{e}^{-}(y)$. Hence the proof is now completed.
\end{proof}
\section{\protect\Large Conclusion}\label{section conclution}
Combining the notions of bipolar fuzzy sets and soft sets leads to the concept of bipolar fuzzy soft sets, that is a more accurate tool for modeling the discussed concepts. From the algebraic point of view, the author \cite{Dehghan Int BF Soft HVS} introduced bipolar fuzzy soft hypervector spaces and studied their basic properties. At the end of the mentioned paper some idea were presented, which some of them were investigated in this paper, such as equivalent conditions of fuzzy bipolar soft hypervector space, bipolar fuzzy soft hypervector spaces generated by a bipolar fuzzy soft set and normal bipolar fuzzy soft hypervector spaces. With this information, the following topics can be studied in the future:

- Cosets of bipolar fuzzy soft hypervector spaces,

- Bipolar fuzzy soft sets on quotient hypervector spaces,

- Finding the applications of the introduced structure in decision making,

- Application of bipolar fuzzy soft sets over other algebraic structures/ hyperstructures.

- Generalized bipolar fuzzy soft hypervector spaces based on bipolar fuzzy points,

- Dimension of bipolar fuzzy soft hypervector spaces.

\vspace{0.7cm}
\noindent {\small\emph{Department of Mathematics, Faculty of Basic Sciences, University of Bojnord, Bojnord, Iran}}

\noindent \emph{dehghan@ub.ac.ir}

\end{document}